\documentclass[leqno]{amsart}
\usepackage{amsmath,amsfonts,amsthm,amssymb,indentfirst,epic,url}

\newcommand{\<}{\kern.0833em}

\newtheorem{theorem}{Theorem}
\newtheorem{lemma}[theorem]{Lemma}
\newtheorem{corollary}[theorem]{Corollary}
\newtheorem{proposition}[theorem]{Proposition}
\newtheorem{definition}[theorem]{Definition}
\newtheorem{question}[theorem]{Question}

\newcommand{\U}{\mathcal{U}}
\newcommand{\F}{\mathcal{F}}

\newcommand{\card}{\mathrm{card}}
\newcommand{\prim}{$\!'\!$}
\newcommand{\idrk}{\mathrm{idp\mbox{-}rk}}
\newcommand{\strt}[1][1.7]{\vrule width0pt height0pt depth#1pt}
\newcommand{\R}{\mathbb{R}}
\newcommand{\C}{\mathbb{C}}

\begin{document}

\begin{center}
\texttt{Comments, corrections, and related references welcomed!}\\[.5em]
{\TeX}ed \today
\vspace{2em}
\end{center}

\title{Homomorphisms on infinite direct product algebras,
especially~Lie~algebras}%
\thanks{After publication of this note, updates, errata, related
references etc., if found, will be recorded at
\url{http://math.berkeley.edu/~gbergman/papers/}.
}

\keywords{homomorphic images of infinite direct products of
nonassociative algebras; simple, solvable, and nilpotent
Lie algebras; maps that factor through ultraproducts;
measurable cardinals%
}

\author{George M. Bergman}
\address[G. Bergman]{University of California\\
Berkeley, CA 94720-3840, USA}
\email{gbergman@math.berkeley.edu}

\author{Nazih Nahlus}
\address[N. Nahlus]{American University of Beirut\\
Beirut, Lebanon}
\email{nahlus@aub.edu.lb}

\begin{abstract}
We study surjective homomorphisms
$f:\prod_I A_i \to B$ of not-necessarily-associative
algebras over a commutative ring $k,$ for $I$ a generally infinite set;
especially when $k$ is a field and
$B$ is countable-dimensional over $k.$

Our results have the
following consequences when $k$ is an infinite field,
the algebras are Lie algebras, and $B$ is finite-dimensional:

If all the Lie algebras $A_i$ are solvable, then so is $B.$

If all the Lie algebras $A_i$ are nilpotent, then so is $B.$

If $k$ is not of characteristic~$2$ or~$3,$ and
all the Lie algebras $A_i$ are finite-dimensional and are direct
products of simple algebras, then (i) so is $B,$
(ii) $f$ splits, and (iii) under a weak cardinality
bound on $I,$ $f$ is continuous in the pro-discrete topology.
A key fact used in getting (i)-(iii) is
that over any such field, every finite-dimensional simple Lie
algebra $L$ can be written $L=[x_1,L]+[x_2,L]$ for
some $x_1,\,x_2\in L,$ which we prove from a
recent result of J.\,M.\,Bois.

The general technique of the paper involves studying
conditions under which a homomorphism on $\prod_I A_i$
must factor through the direct product of finitely many
ultraproducts of the $A_i.$

Several examples are given, and open questions noted.
\end{abstract}
\maketitle

\section{Introduction.}\label{S.intro}
In this note, an {\em algebra} over a commutative associative unital
ring $k$ means a $\!k\!$-module $A$ given with a $\!k\!$-bilinear
multiplication $A\times A\to A,$ which we do not
assume associative or unital.
We shall assume $k$ fixed, and ``algebra'' will
mean ``$\!k\!$-algebra'' unless another base ring is specified.
``Countable'' will be used in the broad
sense, which includes ``finite''.
``Direct product'' will be used in the sense sometimes
called ``complete direct product''.

Let us sketch our method of approach in a somewhat simpler
case than we will eventually be considering.
It is easy to show that if an algebra $B$ is not a nontrivial direct
product, and has no nonzero elements annihilating all of $B,$ then
any surjective homomorphism $f: A_1\times A_2\to B$ from a direct
product of two algebras onto $B$ must factor through the projection
onto one of $A_1$ or $A_2.$
It follows that for such $B,$ a homomorphism
$f: \prod_I A_i\to B$ from an arbitrary direct product onto $B$
will factor through an ultraproduct $\prod_I A_i\,/\,\U.$

For this to be useful, we need to know something
about such ultraproducts.
Assume $k$ a field.
There are three cases:

First, $\U$ may be a principal ultrafilter.
Then $\prod_I A_i\,/\,\U$ can be identified with one
of the $A_i,$ and $f$ factors through the projection to that algebra.

Second, $\U$ may be a nonprincipal ultrafilter that is not
$\!\card(k)^+\!$-complete.
Then $\prod_I A_i\,/\,\U$ is an algebra
over the ultrapower $K = k^I/\,\U,$ and for such $\U,$ $K$ is
an uncountable-dimensional extension field of $k$
(Theorem~\ref{T.bigger}).
We shall see (Proposition~\ref{P.k->K}) that if we map
a $\!K\!$-algebra $A$ onto a $\!k\!$-algebra $B$ having
nonzero multiplication, uncountable dimensionality of $K$
forces $B$ to be uncountable-dimensional over $k$ as well.
Hence, if we restrict attention to maps
onto countable-dimensional $B,$ this case does not occur.

Finally, $\U$ may be a nonprincipal but
$\!\card(k)^+\!$-complete ultrafilter.
If $k$ is finite, $\!\card(k)^+\!$-completeness is vacuous
(implicit in the definition of an ultrafilter), and we cannot
prove much in that case; we mainly note open questions.
If $k$ is infinite, on the other hand, the
$\!\card(k)^+\!$-complete case ``almost'' does not occur:
For such an ultrafilter to exist,
the index set $I$ must be of cardinality at least an
uncountable measurable cardinal, and it is known that if
such cardinals exist, they must be extremely large
(``inaccessible'', and more), and that the
nonexistence of such cardinals is consistent with ZFC, the
standard axiom system of set theory.

Nevertheless, if $\U$ {\em is} such an ultrafilter,
the behavior of $\prod_I A_i\,/\,\U$ is almost as good as
when $\U$ is principal.
This case can be subdivided in two: If the dimensions
of the $A_i$ as $\!k\!$-algebras are not themselves extremely large
cardinals, then that ultraproduct will again be isomorphic
(though not by a projection) to one of the $A_i$
(Theorem~\ref{T.nobigger}), and so will inherit all properties
assumed for these.
Without any restriction on the dimensions of the $A_i,$ the ultraproduct
will still satisfy many important properties that hold on the $A_i,$
e.g., simplicity, nilpotence, or (in the Lie case) solvability
(Propositions~\ref{P.ids} and~\ref{P.simple}),
again allowing us to get strong conclusions about the image~$B$ of $f.$

Fortunately, the proofs of our main results do not require separate
consideration of all these cases, but mainly the distinction
between the $\!\card(k)^+\!$-complete case (which includes the principal
case), and the non-$\!\card(k)^+\!$-complete case
(which, as indicated, will be ruled out under appropriate hypotheses).

The above sketch assumed that $B$ was not a nontrivial
direct product and had no nonzero elements annihilating all of $B.$
We use these hypotheses in the early sections
of the paper, but introduce in \S\ref{S.prod_gen} a weaker
hypothesis on $B$ (``chain condition on almost direct factors'')
yielding more general statements.

In~\S\ref{S.tangents}, we give some tangential results
on concepts introduced in earlier sections.
Some examples showing the need for various of the hypotheses
in our results are collected in~\S\ref{S.examples}.
In~\S\ref{S.questions}, we note some open questions and
directions for further investigation.

Standard definitions and facts about ultrafilters, ultraproducts, and
ultrapowers, assumed from \S\ref{S.via_ultra} on,
are reviewed in an appendix,~\S\ref{S.ultra_defs}.
The more exotic topics of $\!\kappa\!$-complete ultrafilters and
uncountable measurable cardinals
are presented in another appendix,~\S\ref{S.size_ultra},
and used from \S\ref{S.Lie_early} on.

The results proved here about Lie algebras were conjectured several
years ago by the second author (under the assumption that the $A_i$
were finite-dimensional, and the base field $k$ algebraically closed
of characteristic~$0).$

The authors are indebted to Jean-Marie Bois,
Siemion Fajtlowicz, Karl H. Hofmann, Otto Kegel,
A.\,W.\,Knapp, Kamal Makdisi, Donald Passman, Alexander Premet,
and the referree
for many helpful comments and pointers related to this material.

In \cite{prod_Lie2} we obtain
results similar to those of this note -- though stronger in some
ways and weaker in others, and by a rather different approach.

\section{Some preliminaries.}\label{S.prelim}
Our statements in the above sketch, on when a map
$f: A_1\times A_2\to B$ must factor through $A_1$ or $A_2$ on the
one hand, and on when a surjective homomorphism $A\to B$ of
$\!k\!$-algebras, such that $A$ admits a structure of
$\!K\!$-algebra for a ``large'' extension field $K$ of $k,$
forces $B$ to be large, on the other,
both had hypotheses restricting the amount
of ``zero multiplication'' in the structure of $B.$
To see why such limitations are needed, note that if $k$ is a field, and
the $A_i$ are $\!k\!$-algebras whose multiplications are the zero map,
and $B$ is, say, the one-dimensional zero-multiplication
$\!k\!$-algebra, then homomorphisms $\prod_I A_i\to B$
are arbitrary linear functionals on $\prod_I A_i,$ and these need not
satisfy the conclusions of either statement.
Homomorphisms based on zero multiplication are, from our point of
view, very unruly, and we shall ``work around'' that
phenomenon in various ways throughout this paper.

To refer conveniently to that phenomenon, let us make

\begin{definition}\label{D.Z}
If $A$ is an algebra, we define its {\em total annihilator ideal} to be
\begin{equation}\begin{minipage}[c]{35pc}\label{d.Z(A)}
$Z(A)~\,=~\,\{x\in A\mid xA=Ax=\{0\}\}.$
\end{minipage}\end{equation}
\end{definition}

When $A$ is a Lie algebra, $Z(A)$ is thus the {\em center} of $A,$
and our notation agrees with standard notation for the center.
(But when $A$ is an associative algebra, this notation conflicts
with the common notation for the center of $A.)$

Let us also make explicit that the definition of simple
algebra excludes zero multiplication:

\begin{definition}\label{D.simple}
An algebra $A$ will be called {\em simple} if it is nonzero, has
nonzero multiplication, and has no proper nonzero homomorphic images.
\end{definition}

A simple algebra $A$ must be
idempotent, $AA=A,$ and have zero total annihilator, $Z(A)=\{0\},$
since $AA$ and $Z(A)$ are always ideals.

By $AA,$ above, we of course mean the set of sums of products
of pairs of elements of $A.$
More generally,

\begin{definition}\label{D.A'A''}
For any $\!k\!$-submodules $A',\ A''$ of an algebra $A,$ we will
denote by $A'A''$ the $\!k\!$-submodule of $A$ consisting of all sums
of products $a'a''$ $(a'\in A',\,a''\in A'').$
\end{definition}

Below, we will always write $AA$ rather than $A^2,$ to avoid
confusion with $A\times A.$

The next lemma shows that the total annihilator ideal leads to a way
that algebra homomorphisms can be ``perturbed'',
which we will have to take account of in many of our results.
(In this lemma, we explicitly write
``$\!k\!$-algebra homomorphism'' because a $\!k\!$-module
homomorphism is also mentioned.
Elsewhere, ``homomorphism'' will be understood to mean
$\!k\!$-algebra homomorphism unless the contrary is stated.)

\begin{lemma}\label{L.fgh}
Let $A$ and $B$ be $\!k\!$-algebras, $f:A\to B$ a $\!k\!$-algebra
homomorphism, and $h:A\to Z(B)$ a $\!k\!$-module homomorphism.
Then the following conditions are equivalent:

\textup{(i)} $f+h$ is a $\!k\!$-algebra homomorphism.

\textup{(ii)} $AA\ \subseteq\ \ker(h).$

\textup{(iii)} $h$ is a $\!k\!$-algebra homomorphism.
\end{lemma}

\begin{proof}
We will show (i)$\!\iff\!$(ii).
The equivalence (iii)$\!\iff\!$(ii) will follow
as the $f=0$ case of that result.

Since $f+h$ is $\!k\!$-linear,~(i) will hold if and only
if $f+h$ respects multiplication, i.e., if and only if for all
$a,a'\in A$ we have $f(aa')+h(aa')=(f(a)+h(a))(f(a')+h(a')).$
Subtracting from this the equation $f(aa')=f(a)f(a'),$ and noting that
all the terms remaining on the right are two-fold products
in which at least one factor is a value of $h,$ and hence lies
in $Z(B),$ making those products $0,$
we see that~(i) is equivalent to $h(aa')=0,$ i.e.,~(ii).
\end{proof}

Here, in somewhat sharpened form, is the fact stated in
the second paragraph of \S\ref{S.intro}.

\begin{lemma}\label{L.A1A2->B}
Let $B$ be a nonzero algebra.
Then the following conditions are equivalent:

\textup{(i)} \ Every surjective homomorphism
$f: A_1\times A_2\to B$ from a direct product of two algebras
onto $B$ factors through the projection
of $A_1\times A_2$ onto $A_1,$ or through the projection onto $A_2.$

\textup{(ii)} $B$ is not a sum $B_1 + B_2$ of two
nonzero mutually annihilating subalgebras, i.e., nonzero
subalgebras $B_1,\,B_2$ such that $B_1 B_2=B_2 B_1=\{0\}.$

\textup{(iii)} $Z(B)=\{0\},$ and $B$ is not a direct product of
two nonzero subalgebras.

In particular, these conditions hold when $B$ is a simple algebra.
\end{lemma}

\begin{proof}
We shall show that
$\neg\mbox{(i)}\iff\neg\mbox{(ii)}\iff\neg\mbox{(iii)}.$

If a surjective homomorphism $f: A_1\times A_2\to B$
does not factor through either projection
map, then $B_1=f(A_1)$ and $B_2=f(A_2)$ are both nonzero,
and so give a counterexample to~(ii).
Conversely, given $B_1$ and $B_2$ as in $\neg\mbox{(ii)},$
the map $B_1\times B_2\to B$ given by the sum of the inclusions
will be an algebra homomorphism that establishes~$\neg\mbox{(i)}.$

For $B_1$ and $B_2$ as in $\neg\mbox{(ii)}$
we get $\neg\mbox{(iii)}$ by noting that
if they have nonzero intersection, that intersection is a
nonzero submodule of $Z(B),$ while if they have zero
intersection, then $B\cong B_1\times B_2$ as algebras.

Finally, in the situation of $\neg\mbox{(iii)},$ if $Z(B)\neq\{0\}$
then the equation $B=B+Z(B)$ yields $\neg\mbox{(ii)},$
while if $B$ is a direct product of nonzero subalgebras,
then those subalgebras give the required $B_1$ and $B_2.$

The final assertion holds because a simple algebra satisfies~{(iii)}.
\end{proof}

The next lemma strengthens the above result a bit, so as to give
interesting information in the case where $Z(B)\neq\{0\}$ as well.
We shall not use the result in this form, but it will eventually
motivate the transition to the approach of \S\ref{S.prod_gen}.

We remark that the implication (iii\prim)$\!\implies\!$(ii\prim) below
is not a trivial consequence of (iii)$\!\implies\!$(ii) above,
because dividing an algebra $B$ by its total annihilator ideal
$Z(B)$ does not in general produce an algebra with zero total
annihilator ideal, to which we could apply the latter result.

\begin{lemma}\label{L.A1A2->B'}
Let $B$ be a nonzero algebra.
Then the conditions~\textup{(i\prim)} and~\textup{(ii\prim)} below
are equivalent.

\textup{(i\prim)} \ If a homomorphism $f: A_1\times A_2\to B,$
when composed with the natural map $B\to B/Z(B),$ gives a
surjection, then that composite map factors through the projection
of $A_1\times A_2$ onto $A_1,$ or onto $A_2.$

\textup{(ii\prim)} $B$ is not equal to the sum $B_1 + B_2$ of
two mutually annihilating subalgebras neither of which
is contained in $Z(B).$

Moreover, the above conditions are {\em implied} by

\textup{(iii\prim)} $B/Z(B)$ is not a direct product of
two nonzero subalgebras.
\end{lemma}

\begin{proof}[Sketch of proof]
The surjectivity condition in the hypothesis
of~(i\prim) says that $f(A_1)+f(A_2)+Z(B)=B.$
With this in mind, we can get the equivalence of~(i\prim) and~(ii\prim)
as in the preceding result:
given a counterexample to~(i\prim), the subalgebras $B_1=f(A_1)+Z(B)$
and $B_2=f(A_2)+Z(B)$ give a counterexample to~(ii\prim), while
given a counterexample to~(ii\prim), the
induced map $B_1\times B_2\to B$ is a counterexample to~(i\prim).

We complete the proof by showing that
$\neg\mbox{(ii\prim)}\implies\neg\mbox{(iii\prim)}.$
Given $B_1$ and $B_2$ contradicting~(ii\prim),
enlarge them to $B_1+Z(B)$ and $B_2+Z(B)$ if necessary, so that they
each contain $Z(B);$ this clearly preserves the conditions assumed.
It is now easy to show that their images in $B/Z(B)$
contradict~(iii\prim):
Those images are subalgebras which sum to the whole algebra
and annihilate one another on both sides,
so it suffices to show that they have zero intersection.
Since $B_1$ and $B_2$ both contain $Z(B),$ an element of
the intersection of their images in $B/Z(B)$
will arise from an element $x\in B_1\cap B_2.$
But such an element will annihilate both $B_2$ and $B_1,$
hence will annihilate their sum, $B,$ i.e., it will lie in $Z(B),$
so the element of $B/Z(B)$ that it yields is indeed zero.
\end{proof}

The above implication (iii\prim)$\!\implies\!$(ii\prim) is not
reversible.
For example, let $B$ be the associative or Lie
algebra spanned by the matrix units $e_{12},\ e_{13},\ e_{23}$
within the algebra $M_3(k)$ of $3\times 3$ matrices over a field $k.$
(As a Lie algebra,
$B$ is sometimes called the {\em Heisenberg algebra\/}.)
We find that $Z(B)= k\,e_{13},$ and that $B/Z(B)$ is a
$\!2\!$-dimensional $\!k\!$-vector space with zero multiplication, hence
is the direct product of any two one-dimensional subspaces;
so~(iii\prim) fails.
But we claim that~(ii\prim) holds.
Indeed, if $B=B_1+B_2,$ and neither summand is contained in $Z(B),$
then we can find $b_1\in B_1$ and $b_2\in B_2$ which are
linearly independent modulo $Z(B).$
It is then not hard to show that in $[b_1, b_2]=b_1 b_2 - b_2 b_1,$
the coefficient of $e_{13}$ is given by a determinant of
the coefficients of $e_{12}$ and $e_{23}$ in those two elements,
and hence is nonzero;
so $B_1$ and $B_2$ do not annihilate one another, either as
Lie or as associative algebras.

By Lemma~\ref{L.A1A2->B'}, condition~(i\prim) also
holds for this example.
Now a homomorphism from a direct product algebra $A_1\times A_2$
onto the zero-multiplication algebra $B/Z(B)$ can fail to
factor through the projection onto $A_1$ or $A_2$ (e.g., when the latter
are each the one-dimensional zero-multiplication algebra).
Yet condition~(i\prim) shows that {\em if} such a homomorphism arises
as a composite $A_1\times A_2\to B\to B/Z(B),$ it must so factor.

What this example shows us is that though $Z(B)$ is ``trivial'', in
that its elements have zero multiplication with everything, it
cannot be ignored in studying the multiplicative structure of $B$
and the properties of homomorphisms onto $B,$
because elements outside it can have nonzero product lying in it.

\section{Factoring homomorphisms through ultraproducts.}\label{S.via_ultra}

Suppose an algebra $B$ satisfies the equivalent conditions of
Lemma~\ref{L.A1A2->B}, and we map an infinite direct product
$\prod_I A_i$ onto $B.$
Then, since for every subset $J\subseteq I$
we have $\prod_I A_i\cong(\prod_J A_i)\times(\prod_{I-J} A_i),$
Lemma~\ref{L.A1A2->B} gives us a vast family of factorizations
of our homomorphism.
How these fit together is described (in a general set-theoretic
setting) in the next lemma.
In stating it, we assume acquaintance with the concepts of filter,
ultrafilter, reduced product and ultraproduct,
summarized in~\S\ref{S.ultra_defs}.

\begin{lemma}\label{L.via_ultra}
Suppose $(A_i)_{i\in I}$ is a family of nonempty sets,
$B$ is a set, and $f:A=\prod_I A_i\to B$ is a set map,
whose image has more than one element.
Then the following conditions are equivalent:

\textup{(a)} For every subset $J\subseteq I,$ the map $f$
factors either through the projection $A\to\prod_{i\in J}A_i,$
or through the projection $A\to\prod_{i\in I-J}A_i.$

\textup{(b)} The map $f$ factors through the natural map
$A\to A\,/\,\U,$ where $\U$ is an ultrafilter on the index set $I,$
and $A\,/\,\U=\prod_I A_i\,/\,\U$ denotes the ultraproduct of the $A_i$
with respect to this ultrafilter.

When this holds, the ultrafilter $\U$ is uniquely determined by $f.$
\end{lemma}

\begin{proof}
Not yet assuming either~(a) or~(b), but only the initial hypothesis, let
\begin{equation}\begin{minipage}[c]{35pc}\label{d.U=}
$\F\ =\ \{J\subseteq I\mid$
$f$ factors through the projection $A\to\prod_{i\in J}A_i\}.$
\end{minipage}\end{equation}
Thus, a subset $J\subseteq I$ belongs
to $\F$ if and only if for all $a=(a_i)\in A,$
$f(a)$ is unchanged on making arbitrary changes in the coordinates
of $a$ indexed by the elements of the {\em complementary} set $I-J.$
Now if the value of $f(a)$ is unchanged on changing coordinates
lying in a given subset, it is unchanged on changing coordinates
in any smaller subset; and if is unchanged on changing coordinates
in each of two subsets, then it is unchanged on changing coordinates
in the union of those two sets.
Translating these observations into statements about the family $\F$ of
complements of sets with that property, we see that $\F$ is closed
under intersections and enlargement, i.e., $\F$ is a filter on $I.$

Looking at the definition of the reduced product of a family of
sets with respect to a filter on the index set, we see
that $\F$ is the largest filter such that $f$ factors through
the natural map of $A$ to the reduced product $A\,/\,\F.$

The fact that the image of $f$ has more than one element shows
that the value of $f(a)$ is not unchanged under arbitrary modification
of all coordinates of $a;$ so $\F$ does not contain the empty set,
i.e., it is a proper filter.

Finally, we note that condition~(a) is equivalent to saying that for
each $J\subseteq I,$ either $J$ or its complement lies in $\F;$
i.e., that $\F$ is an ultrafilter, which we rename $\U.$
The equivalence of~(a) and~(b), and the final assertion, now follow.
\end{proof}

Combining the above with Lemma~\ref{L.A1A2->B}, we get

\begin{proposition}\label{P.via_ultra}
The equivalent conditions~\textup{(i)-(iii)} of Lemma~\ref{L.A1A2->B}
on a nonzero algebra $B$ are also equivalent~to:

\textup{(iv)} Every surjective homomorphism
$f: \prod_I A_i\to B$ from an arbitrary direct product of algebras
to $B$ factors through the natural map of that product onto the
ultraproduct $\prod_I A_i\,/\,\U,$ for some ultrafilter $\U$ on~$I.$

When this holds, the ultrafilter $\U$ is uniquely determined by $f.$
\end{proposition}

\begin{proof}
Since~(i) is the $I=\{1,2\}$ case of~(iv), we have
(iv)$\!\implies\!$(i).

Conversely, if $B$ satisfies~(i), and we have a homomorphism
$f:A=\prod_I A_i\to B,$ then for every $J\subseteq I$ we can apply~(i)
to the decomposition $A=(\prod_J A_i)\times(\prod_{I-J} A_i),$ and
conclude that $f$ factors through the projection of $A$
to one of these subproducts.
Lemma~\ref{L.via_ultra} now yields~(iv), and the final assertion.
\end{proof}

\section{Extending algebra structures.}\label{S.k->K}

We now come to the other tool referred to in \S\ref{S.intro}.

\begin{proposition}\label{P.k->K}
Suppose $\varphi: k\to K$ is a homomorphism of commutative rings,
$A$ is a $\!K\!$-algebra, $B$ is a $\!k\!$-algebra,
and $f: A\to B$ is a surjective homomorphism as
$\!k\!$-algebras \textup{(}under
the $\!k\!$-algebra structure on $A$ induced by
its $\!K\!$-algebra structure\textup{)}.

Then the kernel of the composite map $A\to B\to B/Z(B)$ is an ideal
of $A,$ not only as a $\!k\!$-algebra, but as a $\!K\!$-algebra.
Hence $B/Z(B)$ acquires a $\!K\!$-algebra structure \textup{(}unique
for the property of making
that composite map a homomorphism of $\!K\!$-algebras\textup{)}.
\end{proposition}

\begin{proof}
It will suffice to show that for $a\in A$ and $c\in K,$
if $f(a)\in Z(B),$ then $f(c\,a)\in Z(B).$
So we must show that $f(c\,a)$ annihilates on both sides an arbitrary
element
of $B,$ which by surjectivity of $f$ we can write $f(a')$ $(a'\in A).$
To do this, we compute: $f(c\,a)f(a')=f(c\,a\,a')=f(a)f(c\,a')=0,$
the last step by the assumption that $f(a)\in Z(B).$
The same calculation works for the product in the opposite order,
completing the proof.
\end{proof}

Remark:
If $A$ were unital, then any {\em ring}\/-theoretic ideal of $A,$
being closed under multiplication by $K\cdot 1,$ would be a
$\!K\!$-algebra ideal.
In that situation, moreover, $Z(B)$ would be trivial, since no
nonzero element of $B$ could be annihilated by $f(1)$ on either side.
The above result shows that somehow, lacking unitality of $A,$ we
can make up for it at the other end by dividing out by~$Z(B).$

\section{First results on homomorphic images of infinite direct product algebras.}\label{S.Lie_early}
In this section we will use the above tools to obtain a couple
of results on homomorphic images of direct products of Lie
and other algebras, under the assumption that
$\card(I)$ is less than any measurable cardinal $>\card(k).$
(That condition holds vacuously, of course,
if no such measurable cardinals exist.)
We assume from here on the material of the final
appendix,~\S\ref{S.size_ultra}, on measurable cardinals, and
$\!\kappa^+\!$-complete and non-$\!\kappa^+\!$-complete ultrafilters.

In the first result below, we restrict the field $k$ so that we
can make use of a theorem of G.\,Brown~\cite{Brown} or its variant in
Bourbaki \cite[Ch.VIII, \S11, Exercise~13(b)]{Bourbaki}, from either of
which it follows that in a finite-dimensional simple Lie
algebra over an {\em algebraically closed field of characteristic~$0,$}
every element is a bracket (not merely a sum of brackets,
as it must be in any simple Lie algebra).
In \S\ref{S.simple_Lie} we will say more precisely what Brown
and Bourbaki prove, then bring in a recent result of
J.\,M.\,Bois which allows us to obtain, with some
more work, a stronger result.

\begin{theorem}\label{T.simple<}
Suppose that $k$ is an algebraically closed field of characteristic~$0,$
that $(A_i)_{i\in I}$ is a family of finite-dimensional
simple Lie algebras over $k,$ that the index set $I$ has
cardinality less than any measurable cardinal $>\card(k),$
and that $f:A=\prod_{i\in I} A_i\to B$ is
a surjective homomorphism to a finite-dimensional Lie algebra $B.$

Then $B$ is semisimple, and $f$ factors
as $\prod_I A_i\to A_{i_1}\times\dots\times A_{i_n}\cong B,$
where the arrow is the projection onto the product of a finite
subfamily of the $A_i.$
\textup{(}In particular, $f$ splits, i.e.,
is right-invertible.\textup{)}
\end{theorem}

\begin{proof}
By the result of Brown and Bourbaki cited above, in each $A_i,$
every element is equal to a single bracket.
Hence the same is true in $A=\prod_I A_i,$ and hence in the homomorphic
image $B$ of $A;$ so in particular, $B$ is idempotent: $B=[B,B].$

Now if $B$ is trivial $\!(0\!$-dimensional), the desired result holds
vacuously with $n=0,$ so assume the contrary.
As a nontrivial idempotent Lie algebra, $B$ must have a homomorphism
onto a simple Lie algebra $C.$
By Proposition~\ref{P.via_ultra}, the composite map
$A\to B\to C$ must factor through the projection of $A=\prod_I A_i$
onto an ultraproduct $A\,/\,\U,$ for some ultrafilter $\U$ on $I.$

By our assumption on the cardinality of $I,$ $\U$ cannot be a
nonprincipal $\!\card(k)^+\!$-complete ultrafilter.
If it were nonprincipal and not $\!\card(k)^+\!$-complete,
then the field $K=k^I/\,\U,$ over which $A\,/\,\U$ is an algebra, would
be uncountable-dimensional over $k$ by Theorem~\ref{T.bigger},
so by Proposition~\ref{P.k->K}, the algebra $C/Z(C)=C$ would
acquire a structure of algebra over $K.$
Hence $C$ would be uncountable-dimensional, contradicting our
hypothesis that $B$ is finite-dimensional.

Hence $\U$ must be a {\em principal} ultrafilter, determined by some
$i_1\in I,$ and what Proposition~\ref{P.via_ultra} then tells us
is that the composite $A\to B\to C$ factors through the projection
$A\to A_{i_1}.$

Now if we write $A=A'\times A_{i_1}$ where
$A'=\prod_{i\in I-\{i_1\}} A_i,$ we see that $B=f(A)$
is the sum of the two mutually annihilating subalgebras
$B'=f(A')$ and $f(A_{i_1}).$
The latter subalgebra, since it does not go to zero under the
map $B\to C,$ and since $A_{i_1}$ is simple, must be an
isomorphic image of $A_{i_1}.$
In particular, it has trivial center.
But the intersection of two mutually annihilating subalgebras
of a Lie algebra must lie in their centers; so the
subalgebras $B'=f(A')$ and $f(A_{i_1})$ have trivial intersection.
Hence $B=B'\times f(A_{i_1}),$
and since $B=[B,B]$ we must have $B'=[B',B'].$

We now repeat the argument with the map $A'\to B'$ in
place of $A\to B.$
By induction on the dimension of $B,$ the ``left-over'' part
(which at this first stage we have called $B')$ must, after finitely
many iterations, become zero, and we get a description of $f,$ up to
isomorphism, as the projection of $A=\prod_I A_i$ onto a finite
subproduct $A_{i_1}\times\dots\times A_{i_n}\cong B.$

Such a projection clearly splits, giving the final assertion.
\end{proof}

Let us next examine homomorphisms on a direct product $\prod_I A_i$
of finite-dimensional {\em solvable} Lie algebras.
We cannot expect that such homomorphisms will in general
factor through finite subproducts, since the solvable Lie algebras
include the abelian ones, which we noted at the beginning of
\S\ref{S.prelim} (under the description ``zero-multiplication
algebras'') can have homomorphisms on their
direct products showing very unruly behavior.
It is nevertheless reasonable to hope that a finite-dimensional
homomorphic image of a direct product of solvable Lie algebras will
be solvable.

Among finite-dimensional Lie algebras, the solvable
ones can be characterized in several ways:
When the base field has characteristic~$0,$ they are those
admitting no homomorphisms onto simple Lie algebras.
In general they are those containing no nonzero idempotent
subalgebras, and also those satisfying one of a certain sequence of
(successively weaker) identities.
Since these conditions do not remain equivalent if one deletes
finite-dimensionality, or the assumption that the
algebras be Lie, or in the finite-dimensional Lie case, that $k$ have
characteristic~$0,$ the statement we hope to
obtain has several possible formulations for general algebras,
all of potential interest.
Of these, the one in terms of nonexistence of homomorphisms
onto simple algebras is ready-made for a proof using
Proposition~\ref{P.via_ultra}.
We shall obtain below a result for general algebras in arbitrary
characteristic based on that proposition, then as a corollary get the
desired statement on finite-dimensional solvable
Lie algebras in characteristic $0.$

Our method will again require a restriction to avoid
complications involving measurable cardinals.
In \S\ref{S.nilp+}, on the other hand, choosing a different condition
that translates to ``solvable'' in the
finite-dimensional Lie case, we will get a result which for
finite-dimensional Lie algebras yields the same conclusion, without
the restrictions on cardinality and characteristic.
For general algebras, however, that result does not subsume
the result of this section; they are independent.

(The technical reason why the present generalization of solvability
will require a cardinality restriction, but the version
in \S\ref{S.nilp+} will not, is that
we have not been able to prove that the property of
admitting no homomorphisms onto
simple algebras is preserved under countably complete ultraproducts,
but we do have the corresponding statement for countable
disjunctions of identities, Proposition~\ref{P.ids}.)

Note that the next result (on not necessarily Lie algebras)
is stronger than the above motivation might lead
one to expect:  the codomain algebra is only required
(in the final sentence) to be countable-dimensional,
rather than finite-dimensional, and still less is assumed about the
dimensionalities of the $A_i.$

\begin{theorem}\label{T.solv}
Suppose $(A_i)_{i\in I}$ is a family of algebras over an infinite
field $k,$ such that no $A_i$ admits a homomorphism onto a
simple algebra \textup{(}or more generally,
such that no $A_i$ admits a homomorphism onto a
countable-dimensional simple algebra\textup{)}.

Assume further that, if there are measurable cardinals
greater than $\card(k),$ then either $\card(I)$ or the supremum of the
dimensions of all the $A_i$ is less than all such cardinals.

Then $\prod_I A_i$ admits no homomorphism onto a
countable-dimensional simple \mbox{algebra}.
\end{theorem}

\begin{proof}
Suppose, by way of contradiction, that $f: \prod_I A_i\to B$ is a
homomorphism onto a countable-dimensional simple algebra.
Lemma~\ref{L.A1A2->B} and Proposition~\ref{P.via_ultra}
tell us that for some ultrafilter $\U$ on $I,$ this $f$ factors
through the natural map $\prod_I A_i\to\prod_I A_i\,/\,\U.$

If $\U$ is principal (in which case it is $\!\mu\!$-complete
for every cardinal $\mu),$ or is nonprincipal but
$\!\card(k)^+\!$-complete (in which case, $\card(I)$ must be
greater than or equal to some measurable cardinal $\mu>\card(k),$
so that the bound on the dimensions of the $A_i$ in the second
paragraph of the theorem applies), then Theorem~\ref{T.nobigger} shows
that $\prod_I A_i\,/\,\U$ is isomorphic to one of the
$A_i;$ but by hypothesis, no $A_i$ admits a homomorphism onto
$B,$ a contradiction.

On the other hand, if $\U$ is not $\!\card(k)^+\!$-complete,
we can argue exactly as in the proof of Theorem~\ref{T.simple<},
and conclude that $B$ is uncountable-dimensional, though we
assumed the contrary.

So there exists no such $f,$ as was to be proved.
\end{proof}

Here is the resulting statement about solvable Lie algebras.
The dimensions of the $A_i$ are still almost unrestricted,
but we must make $B$ finite-dimensional to turn the
``no simple images'' condition into solvability.

\begin{corollary}\label{C.solv}
Suppose $(A_i)_{i\in I}$ is a family of solvable Lie algebras
over a field $k$ of characteristic~$0,$ and suppose that,
if there exists a measurable cardinal greater than $\card(k),$
then either $\card(I)$ or the supremum of the dimensions of
the $A_i$ is less than every such cardinal.
\textup{(}E.g., this is automatic if all $A_i$ are
of dimension $\leq$ the continuum; in particular, if they
are finite-dimensional.\textup{)}

Then any finite-dimensional homomorphic image $B$ of
$\prod_I A_i$ is solvable.
\end{corollary}

\begin{proof}
This follows from the preceding theorem, since if a Lie algebra is
solvable, it admits no homomorphism onto a simple Lie algebra, and
the converse holds in the finite-dimensional case in characteristic~$0.$
\end{proof}

\section{Ultraproducts and almost direct factors.}\label{S.prod_gen}

The arguments of the preceding section were based on reducing
the results to be proved to the consideration of homomorphisms from
our infinite product onto simple algebras, and
applying Proposition~\ref{P.via_ultra}.
But there are situations where that method is not enough.
If the $A_i$ are simple non-Lie algebras,
we do not have Brown's theorem available
to tell us that $\prod_I A_i$ is idempotent, from which we
deduced that $B$ had to have a simple homomorphic image.
And if, in our consideration of solvable Lie algebras, we either replace
solvability by nilpotence, or look at characterizations of solvability
applicable in arbitrary characteristic,
we again can't use that argument.
For such purposes, we would like to have some variant of
Proposition~\ref{P.via_ultra} not
burdened with the condition $Z(B)=\{0\}$
of Lemma~\ref{L.A1A2->B}(iii); something more in the
spirit of Lemma~\ref{L.A1A2->B'} than of Lemma~\ref{L.A1A2->B}.

It would also be nice to replace the assumption
of finite-dimensionality of $B$ in the conclusions of
both Theorem~\ref{T.simple<} and Corollary~\ref{C.solv}
by a more general condition.

We develop below a refinement of Proposition~\ref{P.via_ultra} in
line with these two ideas.
We will use the following concept, motivated by
condition (ii\prim) of Lemma~\ref{L.A1A2->B'}.

\begin{definition}\label{D.almost}
For any algebra $A,$ an {\em almost direct decomposition} of $A$
will mean an expression $A=\nolinebreak B+B',$ where $B,\ B'$ are
ideals of $A,$ and each is the two-sided annihilator of the other.
In this situation, $B$ and $B'$ will be called
{\em almost direct factors} of $A,$ and each will be called
the {\em complementary} factor to the other.
\end{definition}

Remarks: Any algebra $A$ has a smallest almost direct factor, $Z(A);$
its complement, the largest almost direct factor, is $A.$

If $A=B_1+B_2$ is an almost direct decomposition, then
$Z(B_1)=Z(B_2)=Z(A).$
Indeed, because $B_1$ is the annihilator of $B_2$ we
have $Z(A)\subseteq B_1,$ and hence
$Z(A)\subseteq Z(B_1);$ conversely, any $a\in Z(B_1)$
annihilates both $B_1$ and $B_2,$ hence
annihilates $A,$ i.e., lies in $Z(A).$

If $Z(A)=\{0\},$ the almost direct decompositions of $A$ are its
(internal) pairwise direct product decompositions as an algebra.
If $Z(A)\neq\{0\},$ an almost direct
decomposition is never a direct product decomposition, since the two
factors in the decomposition intersect in $Z(A).$
However, an almost direct decomposition of $A$
does induce a direct product decomposition of $A/Z(A),$ as shown in
the proof of Lemma~\ref{L.A1A2->B'}(iii\prim)$\!\implies\!$(ii\prim).
On the other hand, not every direct product decomposition of $A/Z(A)$
need arise in that way, as shown by the $3\times 3$ matrix example
following that lemma.

An almost direct factor of an almost direct factor of an algebra $A$ is
easily seen to be, itself, an almost direct factor of $A.$
If we perform finitely many such
successive almost direct decompositions, we get a decomposition
$A=B_1+\dots+B_n$ as a sum of ideals each of which is the annihilator
of the sum of the rest.
We may call such an expression an almost direct
decomposition into several almost direct factors.

The importance for us of almost direct decompositions lies in

\begin{lemma}\label{L.prod->>}
If $f: A_1\times A_2\to B$ is a surjective homomorphism of algebras,
or more generally, a homomorphism satisfying
$f(A_1\times A_2)+Z(B)=B,$ then $f(A_1)+Z(B)$ and $f(A_2)+Z(B)$ are
complementary almost direct factors of $B.$
\end{lemma}

\begin{proof}
By assumption, $f(A_1)+Z(B)$ and $f(A_2)+Z(B)$ sum to $B,$ so
it remains to prove that they are ideals of $B,$ and are mutual
two-sided annihilators.
By symmetry, it suffices to show that $f(A_1)+Z(B)$ is an ideal,
and is the two-sided annihilator of $f(A_2)+Z(B).$

Since $A_1$ is an ideal of $A_1\times A_2,$ its
image $f(A_1)$ is an ideal of
$f(A_1\times A_2)+Z(B)=B,$ hence so is $f(A_1)+Z(B).$
Since $A_1$ and $A_2$ annihilate one another in $A_1\times A_2,$
$f(A_1)+Z(B)$ and $f(A_2)+Z(B)$ annihilate one another in $B.$
So it remains to show that any $b\in B$ which annihilates
$f(A_2)+Z(B)$ lies in $f(A_1)+Z(B).$

Let us write such an element $b$ as $f(a_1)+f(a_2)+z,$ with
$a_i\in A_i,\,z\in Z(B).$
Since $f(a_1)$ and $z$ automatically annihilate $f(A_2)+Z(B),$
and by assumption $b=f(a_1)+f(a_2)+z$ does,
it follows that $f(a_2)$ does.
But as a member of $f(A_2),$ it also annihilates $f(A_1),$ hence it
annihilates all of $f(A_1)+f(A_2)+Z(B)=B,$ i.e., lies in $Z(B).$
Hence the expression $b=f(a_1)+(f(a_2)+z)$ expresses $b$
as a member of $f(A_1)+Z(B),$ as required.
\end{proof}

Here is the generalization of finite-dimensionality
that we indicated would be helpful in strengthening our results.

\begin{definition}\label{D.CC}
We shall say that an algebra $B$ has
{\em chain condition on almost direct factors}
if it has no infinite strictly ascending chain
$B_1\subsetneq B_2 \subsetneq\dots$ of almost direct factors;
equivalently, if it has no infinite strictly descending chain
$B_1'\supsetneq B_2' \supsetneq\dots$
of almost direct factors; equivalently, if it has
no infinite chain \textup{(}totally ordered set\textup{)}
of almost direct factors.
\end{definition}

We do not know much general theory regarding the above chain condition.
Clearly, it will hold in a ring with chain condition on two-sided
ideals, which is already much weaker than finite-dimensionality.
After the main results of this paper, a couple of results on the
condition will be proved in \S\S\ref{S2.Boole}-\ref{S2.weaker}.
An example in \S\ref{S2.CC} will show that not every finitely
generated associative algebra over a field satisfies it.
We do not know whether every finitely generated Lie algebra over a
field does

Here, now, is our modified version of Proposition~\ref{P.via_ultra}:

\begin{proposition}\label{P.via_ultras/Z}
Suppose $f: \prod_I A_i\to B$ is a surjective homomorphism of
algebras, or more generally, a homomorphism such
that $f(\prod_I A_i)+Z(B)=B;$
and suppose $B$ has chain condition on almost direct factors.
Let us abbreviate $\prod_I A_i$ to $A,$ and write $\pi: B\to B/Z(B)$
for the canonical factor map.

Then there exists a finite family of
distinct ultrafilters $\U_1,\dots,\U_n$
on $I$ such that, if we write $\varphi_m$ for the natural homomorphism
$A\to A\,/\,\U_m$ $(m=1,\dots,n),$ then
the composite map $\pi f: A\to B/Z(B)$ factors through the
map $(\varphi_1,\dots,\varphi_n):
A\to A\,/\,\U_1\times\dots\times A\,/\,\U_n.$
\end{proposition}

\begin{proof}
If $B=Z(B)$ this is vacuous, so assume the contrary.

For any partition $I=J\cup(I-J),$ we have the direct
product decomposition $A=\prod_J A_i\times \prod_{I-J} A_i;$ so
by Lemma~\ref{L.prod->>}, $f(\prod_J A_i)+Z(B)$ is an
almost direct factor of~$B,$
with complementary almost direct factor $f(\prod_{I-J} A_i)+Z(B).$
Inclusions of subsets $J$ give inclusions of
almost direct factors $f(\prod_J A_i)+Z(B),$
so by our assumption of chain condition
on such factors, there must exist some $J_1\subseteq I$
which yields an almost direct factor $B_1$ strictly larger than
$Z(B),$ but such that every subset $J\subseteq J_1$ yields either
the same almost direct factor, $B_1,$ or the trivial
almost direct factor, $Z(B).$

Now by Lemma~\ref{L.prod->>}, for every $J\subseteq J_1$
the ideals $f(\prod_{J} A_i)+Z(B)$ and $f(\prod_{J_1-J}A_i)+Z(B)$
are complementary almost direct factors of $B_1;$
but by our choice of $J_1,$ these
can only be $B_1$ and $Z(B)$ in one or the other order.
We claim that the set
\begin{equation}\begin{minipage}[c]{35pc}\label{d.U(0)}
$\U_1^{(0)}\ =\ \{J\subseteq J_1\mid f(\prod_{J} A_i)+Z(B)=B_1\}$
\end{minipage}\end{equation}
is an ultrafilter on $J_1.$
Indeed, the class of subsets of $J_1$ with the reverse
property, $f(\prod_{J} A_i)+Z(B)=Z(B),$ is closed under pairwise unions
and passage to subsets, so $\U_1^{(0)}$ is closed under
pairwise intersections and enlargements, i.e., is a filter.
Clearly $\emptyset\notin\U_1^{(0)},$ so this filter is proper.
And we have noted that for every $J\subseteq J_1,$ one
of $J$ or $J_1-J$ is in $\U_1^{(0)},$ so it is an ultrafilter.

Having found an ultrafilter that roughly describes the behavior of $f$
as a map from $\prod_{J_1} A_i$ into the almost
direct factor $B_1=f(\prod_{J_1} A_i)+Z(B)$ of $B,$
we now look at the complementary
factor $\prod_{I-J_1} A_i$ of $A,$ which $f$ maps into the complementary
almost direct factor $f(\prod_{I-J_1} A_i)+Z(B)$ of $B.$
If the latter is not $Z(B),$ we can repeat the above process
with this map.
(It was to make this work that we put the ``{\em or more generally}''
clause
into the first sentence of this proposition and the preceding lemma.
If $f$ was assumed surjective to $B,$ this would not guarantee
that its restriction to $\prod_{I-J_1} A_i$ would be surjective to
$f(\prod_{I-J_1} A_i)+Z(B).)$
Thus we get a subset $J_2\subseteq I-J_1$
such that $f(\prod_{J_2} A_i)+Z(B)$ is a minimal nontrivial
almost direct factor of $f(\prod_{I-J_1} A_i)+Z(B),$ and
an ultrafilter $\U_2^{(0)}$ on that subset such that
every member of $\U_2^{(0)}$ induces that same almost direct factor.

Iterating this process, we get a strictly decreasing sequence
of almost direct factors of $B$ associated with the sets
$I,$ $I-J_1,$ $I-J_1-J_2,$ $\dots\,;$ so our chain condition insures
that this iteration cannot continue indefinitely.
Thus, at some stage, say the $\!n\!$-th, our
complementary almost direct factor must be $Z(B),$ so
the factor whose complement it is must be the whole algebra
we are considering at that stage.
Thus, without loss of generality we may, at that stage,
take $J_n=I-J_1-\dots-J_{n-1}$ (rather than some proper subset thereof),
giving us a partition $I=\bigcup_{m=1,\dots,n}J_m;$
and we see that the ideals $f(\prod_{J_m} A_i)+Z(B)$ $(m=1,\dots,n)$
constitute an almost direct decomposition of~$B.$

Now for $m=1,\dots,n,$ let $\U_m$ be the ultrafilter
on $I$ induced by the ultrafilter $\U_m^{(0)}$ on $J_m,$
i.e., $\U_m=\{J\subseteq I\mid J\cap J_m\in\U_m^{(0)}\}.$
For each such $m$ we define a homomorphism
\begin{equation}\begin{minipage}[c]{35pc}\label{d.gm}
$g_m:\ A\,/\,\U_m\ \to\ B/Z(B)$
\end{minipage}\end{equation}
as follows.
Any element of $A\,/\,\U_m$ is the image of some
$a=(a_i)_{i\in I}\in A = \prod_I A_i.$
Let us map this by restriction to $\prod_{J_m} A_i,$ and then by
inclusion into $A;$ this means replacing the components $a_i$ at
indices $i\notin J_m$ by zero, while keeping the components at
indices $i\in J_m$ unchanged.
Map the resulting element by $f$ into $B,$ and then by $\pi$
into~$B/Z(B).$

If we had chosen a different representative $a'=(a'_i)_{i\in I}\in A$
of our element of $A\,/\,\U_m,$ then after restriction to $J_m,$
this would have differed from $(a_i)_{i\in I}$ only
on a subset of $J_m$ that is not in $\U_m^{(0)}.$
But by our construction of $\U_m^{(0)},$
elements with support in such a subset of $J_m$ are mapped into $Z(B)$
by $f;$ so the image under $\pi f$ of the element obtained from $a'$
equals the image under $\pi f$ of the element obtained from $a,$
showing that we have described a well-defined map~(\ref{d.gm}).

It is now routine to verify that $g_m$ is a homomorphism
$A\,/\,\U_m\to B/Z(B),$ with image in $(f(\prod_{J_m} A_i)+Z(B))/Z(B),$
so that the images of $g_1\varphi_1,\ \dots,\ g_n\varphi_n$ annihilate
one another; and that $\pi f=g_1\varphi_1+\dots+g_n\varphi_n,$
so that this map indeed factors through $(\varphi_1,\dots,\varphi_n).$
\end{proof}

The above result says that under the indicated hypotheses we can, in
a certain sense, approximate $f:A\to B$ ``modulo $Z(B)$''
by a homomorphism
that factors through $A\,/\,\U_1\times\dots\times A\,/\,\U_n.$
It is natural to ask whether we can do so in a stronger sense, namely
whether we can express $f$ as a ``perturbation'', of the sort described
by Lemma~\ref{L.fgh}, of a genuine homomorphism $f_1:A\to B$
factoring through $A\,/\,\U_1\times\dots\times A\,/\,\U_n.$

We can do this easily if the ultrafilters $\U_m$ are principal:
if each $\U_m$ is the principal ultrafilter determined
by $i_m\in I,$ one finds that the desired $f_1:A\to B$ can be
obtained by projecting $A$ to $A_{i_1}\times\dots\times A_{i_n}$
regarded as a subalgebra of $A,$ and then mapping by $f$ into $B.$

For nonprincipal $\U_m,$ we do not know whether such a factorization is
always possible; but we shall show that it is whenever $k$ is a field.
Note that what we want is to perturb the given
homomorphism $f$ to a homomorphism $f_1$ whose kernel
contains the kernel of the natural surjection
$(\varphi_1,\dots,\varphi_n):
A\to A\,/\,\U_1\times\dots\times A\,/\,\U_n.$
The following is a general result on when a homomorphism of
algebras over a field has
a perturbation whose kernel contains a prescribed ideal.

\begin{lemma}\label{L.AAcapC}
Let $k$ be a field, $f:A\to B$ any homomorphism
of $\!k\!$-algebras, and $C$ an ideal of $A.$
Then the following conditions are equivalent:\vspace{.2em}

\textup{(i)} There exists a homomorphism $f_1:A\to B$ having
$C$ in its kernel, such that $f-f_1$ is $\!Z(B)\!$-valued.\vspace{.2em}

\textup{(ii)} $f(C)\subseteq Z(B),$ and
$f(AA\,\cap\,C)=\{0\}.$\vspace{.2em}

Moreover, if $B=f(A)+Z(B),$ then the first condition of \textup{(ii)}
is implied by the second.
\end{lemma}

\begin{proof}
To get (i)$\!\implies\!$(ii) (which does not require the assumption
that $k$ is a field), suppose we have an $f_1$ as in~(i).
Then $C$ is carried into $Z(B)$ by $f-f_1,$ but annihilated
by $f_1,$ hence it must be carried into $Z(B)$ by $f,$ giving the
first assertion of~(ii).
Further, Lemma~\ref{L.fgh} (with $h=f_1-f)$
tells us that $f_1-f$ annihilates $AA;$ and by assumption, $f_1$
annihilates $C,$ so $f=f_1-(f_1-f)$ must annihilate $AA\,\cap\,C,$
giving the second assertion.

Conversely, assuming~(ii), the second assertion thereof shows that
the zero map and $-f$ agree on $AA\cap C,$ whence
there exists a unique $\!k\!$-linear map $h:AA+C\to B$ that
agrees with the zero map on $AA$ and with $-f$ on $C;$ and by
the first condition of (ii), it will be $\!Z(B)\!$-valued.
If $k$ is a field, we can extend this as a vector space map
(in an arbitrary way) to a map $h:A\to Z(B).$
Since $h$ annihilates $AA,$ Lemma~\ref{L.fgh} tells us that
$f_1=f+h$ is a $\!k\!$-algebra homomorphism; and since $h$
agrees with $-f$ on $C,$ $f+h$ has $C$ in its kernel.

For the final statement (which again does not require that
$k$ be a field), note that to show that $f(C)\subseteq Z(B)$
is to show that $f(C)$ is annihilated on each side by $B,$
which, if $B=f(A)+Z(B),$ is equivalent to being annihilated
on each side by $f(A).$
But multiplication on either side by $f(A)$ carries $f(C)$ into
$f(AC)+f(CA)\subseteq f(AA\cap C),$ which
is zero assuming the second condition of~(ii).
\end{proof}

\begin{corollary}\label{C.via_ultras+}
Under the hypotheses of Proposition~\ref{P.via_ultras/Z}, if $k$
is a field, then $f$ can be written $f_1+f_0,$ where $f_1$
is a homomorphism $A\to B$ factoring through
$(\varphi_1,\dots,\varphi_n):
A\to A\,/\,\U_1\times\dots\times A\,/\,\U_n,$
and $f_0$ is a homomorphism $A\to Z(B)$ \textup{(}necessarily
factoring through the natural map $A\to A/AA).$
\end{corollary}

\begin{proof}
It is easy to see that $(\varphi_1,\dots,\varphi_n)$
maps $A$ surjectively to $A\,/\,\U_1\times\dots\times A\,/\,\U_n;$
hence a homomorphism on $A$ can be factored through
that map if and only if its kernel contains
$\ker(\varphi_1)\cap\,\dots\,\cap\ker(\varphi_n).$
Hence, by Lemma~\ref{L.AAcapC} (including the final sentence),
with $C$ taken to be that intersection of kernels,
it suffices to show that any $a$ belonging both to
$\ker(\varphi_1)\cap\dots\cap\ker(\varphi_n)$ and
to $AA$ is in $\ker(f).$

Given such an $a,$ let $J=\{i\in I\mid a_i=0\}.$
Thus, letting $A_1=\prod_{J} A_i$ and $A_2=\prod_{I-J} A_i,$
we have $A=A_1\times A_2,$
and $a$ has zero component in the first factor.
Hence since $a\in AA=A_1A_1+A_2A_2,$ we have $a\in A_2A_2.$

Also, since $a$ lies in the kernels of all $\varphi_m,$ its
support $I-J$ belongs to none of the $\U_m.$
Hence $A_2,$ which is also supported on that set, is
likewise contained in the kernels of all the maps $\varphi_m;$ so
by the conclusion of Proposition~\ref{P.via_ultras/Z},
$A_2\subseteq\ker(\pi\,f).$
This says that $f(A_2)\subseteq Z(B);$
hence $f(a)\in f(A_2 A_2)\subseteq Z(B)Z(B)=\{0\},$ as required.

The final parenthetical statement is a case of the general
observation that any homomorphism from an algebra
$A$ to a zero-multiplication algebra factors through $A/AA.$
\end{proof}

If we now add to Corollary~\ref{C.via_ultras+} the
assumptions that the field $k$ is infinite and the
algebra $B$ countable-dimensional, we can again (as
in the proof of Theorem~\ref{T.simple<}) use
Proposition~\ref{P.k->K} and Theorem~\ref{T.bigger} to exclude
the case where any of the ultrafilters $\U_m$ are
non-$\!\card(k)^+\!$-complete, getting

\begin{theorem}\label{T.via_A_i}
Suppose $f:A=\prod_I A_i\to B$ is a surjective homomorphism of algebras
over an infinite field $k,$ where $B$ is countable-dimensional
over $k$ and has chain condition on almost direct factors.

Then there exist finitely many
distinct {\em $\!\card(k)^+\!$-complete} ultrafilters
$\U_1,\dots,\U_n$ on $I$ such that, writing $\varphi_m$ for the natural
homomorphism $A\to A\,/\,\U_m$ $(m=1,\dots,n),$
$f$ can be written $f_1+f_0,$ where
$f_1$ factors through the map $(\varphi_1,\dots,\varphi_n):
A\to A\,/\,\U_1\times\dots\times A\,/\,\U_n,$ and
$f_0$ is a homomorphism $A\to Z(B).$

In particular, if $\card(I)$ is not $\geq$ any
measurable cardinal $>\card(k),$ then each of the $\U_m$
is a principal ultrafilter, so $f_1$ factors through the projection
of $A$ onto the product of finitely many of the $A_i.$

If it is merely assumed that none of the dimensions
$\mathrm{dim}_k(A_i)$ is $\geq$ a measurable cardinal $>\card(k),$
then the algebras
$A\,/\,\U_m$ are, at least, each isomorphic to one of the $A_i.$
\qed
\end{theorem}

(Remark: if $k$ is uncountable, the proof of Theorem~\ref{T.bigger}
shows that for a non-$\!\card(k)^+\!$-complete ultrafilter $\U,$
the field $k^I/\,\U$ will have dimension at least $\card(k)$
over $k.$
So in that case, one can get the above result not only for $B$
countable-dimensional, but for $B$ of any dimensionality $<\card(k).)$

\section{Solvable algebras (version 2), and nilpotent algebras.}\label{S.nilp+}

We can now, as promised, prove a result which,
for the finite-dimensional Lie case, gives essentially the same
conclusion about homomorphic images of direct products of solvable Lie
algebras as Corollary~\ref{C.solv}, but without the
restrictions on $\card(I)$ and $\mathrm{char}(k),$
while for general algebras, it is independent of that result.
We shall also obtain the analogous result for nilpotent algebras.

Our conditions on general algebras will use the following analogs
of the derived series and the lower central series of a Lie algebra.
(We modify slightly a common notation for the latter, to avoid confusion
with the subscripts indexing the factors in our direct products.)

\begin{definition}\label{D.series}
In any algebra $A,$ we define, recursively, $\!k\!$-submodules
$A^{(d)}$ $(d=0,1,\dots)$ and $A_{[d]}$ $(d=1,2,\dots)$ by
\begin{equation}\begin{minipage}[c]{35pc}\label{d.der_series}
$A^{(0)}=A,\qquad A^{(d+1)}= A^{(d)}\,A^{(d)},$
\end{minipage}\end{equation}
\begin{equation}\begin{minipage}[c]{35pc}\label{d.low_cen_series}
$A_{[1]}\,=A,\qquad A_{[d+1]}= A\,A_{[d]}+A_{[d]}A.$
\end{minipage}\end{equation}
We will call $A$ {\em solvable} if $A^{(d)}=\{0\}$ for some $d,$
and {\em nilpotent} if $A_{[d]}=\{0\}$ for some $d.$
\end{definition}

(The concept of nilpotence of a general algebra is standard.
That of solvability is less so, but it appears in \cite[p.17]{Schafer}.
It is not hard to show that the submodules $A_{[d]}$ are ideals,
and that the $A^{(d)}$ are subalgebras, and in the Lie case are
ideals as well; but we shall not need these facts here.)

\begin{theorem}\label{T.solv2}
Suppose $k$ is an infinite field,
and $(A_i)_{i\in I}$ is a family of solvable
$\!k\!$-algebras, in the sense of Definition~\ref{D.series}
\textup{(}e.g., solvable Lie algebras in the
standard sense\textup{)}.
Then any finite-dimensional homomorphic image of
$\prod_I A_i$ is solvable.
\end{theorem}

\begin{proof}
Say $f: A=\prod_I A_i\to B$ is a homomorphism onto a finite-dimensional
algebra.
Since finite-dimensionality implies chain condition
on almost direct factors, Theorem~\ref{T.via_A_i} shows that $B$ is
a sum of finitely many mutually annihilating homomorphic images of
algebras $A\,/\,\U_m,$ where the $\U_m$ are $\!\card(k)^+\!$-complete
ultrafilters on $I,$ together with a subspace $f_0(A)\subseteq Z(B).$

The $\U_m$ are, in particular, countably complete, and
solvability is equivalent to the condition that
an algebra satisfy one of the countable family of identities,
\begin{equation}\begin{minipage}[c]{35pc}\label{d.solvids}
$x=0,\quad x_0\,x_1=0,\quad (x_{00}\,x_{01})\,(x_{10}\,x_{11})=0,\quad
\dots\ .$
\end{minipage}\end{equation}
Hence by Proposition~\ref{P.ids}, the condition of solvability
on the $A_i$ carries over to the algebras $A\,/\,\U_m.$
The subalgebra $f_0(A)\subseteq Z(B)$ clearly also satisfies
the second identity of~(\ref{d.solvids}).
Hence, as the identities of~(\ref{d.solvids}) are successively weaker,
at least one of them will be satisfied by all of the finitely
many algebras $A\,/\,\U_m$ and by~$f_0(A).$

A sum of finitely many mutually annihilating algebras satisfying
a common identity also satisfies that identity, yielding
the asserted solvability.
\end{proof}

Exactly the same method yields

\begin{theorem}\label{T.nilp}
Suppose $(A_i)_{i\in I}$ is a family of nilpotent
algebras \textup{(}e.g., nilpotent Lie algebras\textup{)} over an
infinite field $k.$
Then any finite-dimensional homomorphic image of
$\prod_I A_i$ is nilpotent.\qed
\end{theorem}

In this case, however, a stronger result will be proved in
\cite{pro-np}, by different methods, with ``direct product'' generalized
to ``inverse limit'', and no requirement that $k$ be infinite.
(From that result of \cite{pro-np}, an analog of Theorem~\ref{T.solv2}
for inverse limits of solvable algebras is also deduced,
but only for finite-dimensional Lie algebras
$A_i$ over a field of characteristic~$0,$ for which there is
an easy characterization of solvability in terms of nilpotence.)

\section{Simple algebras -- general results.}\label{S.simple}

We would now like to use Theorem~\ref{T.via_A_i} to get a result
on homomorphic images of products of finite-dimensional {\em simple}
Lie algebras stronger than our earlier Theorem~\ref{T.simple<}.
Simplicity is not, like solvability or nilpotence,
equivalent to a disjunction of identities, but that is not a problem:
like countable disjunctions of identities, it is preserved by
countably complete ultraproducts (Proposition~\ref{P.simple}).
A more serious difficulty is that the preceding proofs
used the fact that $f_0(A),$ a zero-multiplication algebra,
was automatically nilpotent and solvable;
but if $f_0(A)\neq\{0\},$ it will certainly
{\em not} be a product of simple algebras.

By the final observation of Corollary~\ref{C.via_ultras+}, the map $f_0$
of Theorem~\ref{T.via_A_i} can be nontrivial only if $AA\neq A,$
i.e., if $A$ is not idempotent.
Simple algebras $A_i$ {\em are} idempotent.
Does this property carry over to direct products?

To answer this, let us define, for every idempotent
algebra $A,$ its {\em idempotence rank,} $\idrk(A),$ to be
the supremum, over all $a\in A,$ of the least number $m$ of
summands in expressions for $a$ as a sum of products:
\begin{equation}\begin{minipage}[c]{35pc}\label{d.*Sab}
$\idrk(A)\ =\ \sup_{a\in A}(\inf\,\{m\geq 0\mid
(\exists\,b^{(1)},\dots,b^{(m)},\,c^{(1)},\dots,c^{(m)}\in A)
\ \ a=\sum_{h=1}^m b^{(h)} c^{(h)}\}).$
\end{minipage}\end{equation}
This will be a nonnegative integer (positive if $A\neq\{0\}),$
or $\infty;$ it measures the difficulty in asserting
the idempotence of $A$ in a uniform way.
We can now state and prove

\begin{lemma}\label{L.idpt_rk}
For a family of algebras $A_i$ $(i\in I),$
the following conditions are equivalent.\vspace{.2em}

\textup{(i)} $\prod_I A_i$ is idempotent.\vspace{.2em}

\textup{(ii)} Every $A_i$ is idempotent, and there is a natural
number $n$ such that for all but finitely many $i\in I,$
$\idrk(A_i)\leq n.$\vspace{.2em}

When the above equivalent conditions hold,
$\idrk(\prod_I A_i)=\sup_{i\in I} \idrk(A_i).$
\end{lemma}

\begin{proof}
We shall prove (ii)$\!\implies\!$(i) and
$\neg\mbox{(ii)}\implies\neg\mbox{(i)}.$

Given $n$ as in~(ii), consider any $(a_i)_{i\in I}\in A.$
For those $i$ such that $\idrk(A_i)\leq n,$ take a
representation of $a_i$ as a sum of $n$ products,
while for each of the remaining finitely many indices $i,$ take
{\em some} representation of $a_i$ as a sum of products.
(Some of the $A_i$ may have infinite idempotence rank; but they are all
assumed idempotent, so each element $a_i$ can be so written.)
There will be a common upper bound $N$ for the number of summands
in all these representations, yielding a representation for
$(a_i)$ as a sum of $N$ products, proving~(i).

Assuming $\neg\mbox{(ii)},$ note that if not
all $A_i$ are idempotent, then $A$ cannot be.
If they all are idempotent, but there is no finite $n$ bounding
the idempotence ranks of all but finitely many of them,
then it is easy to construct
an $(a_i)\in A$ such that the number of products required to
express the component $a_i$ is unbounded as a function of $i,$ and
to deduce that $(a_i)$ cannot be written as a finite sum of products,
proving $\neg\mbox{(i)}.$

The verification of the final assertion (which we won't
use) is straightforward; one breaks it into two cases,
the case where the idempotence ranks of all
the $A_i$ are finite, so that~(ii) implies that they have a common
finite bound, and the case where at least one is $\infty.$
\end{proof}

Applying this to homomorphic images
of direct products of simple algebras, we can now prove

\begin{theorem}\label{T.simple}
Suppose $k$ is an infinite field, and
$f:A=\prod_I A_i\to B$ is a surjective homomorphism from a
direct product of simple algebras
to a countable-dimensional algebra $B$
having chain condition on almost direct factors.
Then

\textup{(a)} If there is a finite upper bound on all but finitely
many of the values $\idrk(A_i)$ $(i\in I),$ then $B$
is isomorphic to a direct product of finitely many of the $A_i.$

\textup{(b)} Without the assumption of such a bound,
$B$ will be isomorphic to a direct product of finitely many of the
$A_i,$ and one $\!k\!$-vector-space with zero multiplication.

In the situation of \textup{(a)}, the homomorphism $f$
splits \textup{(}has a right inverse\textup{)}.
In the situation of \textup{(b)}, the
composite homomorphism $A\to B\to B/Z(B)$ splits.
\end{theorem}

\begin{proof}
Let $f$ be expressed as in Theorem~\ref{T.via_A_i}.
By Proposition~\ref{P.simple}, all of the $A\,/\,\U_m$
in that description are simple.
A homomorphic image of a finite direct product of simple algebras is
the direct product of some subset of these; so possibly dropping
some of the $A\,/\,\U_m,$ we may assume that the
image of $f_1$ in $B$ is an isomorphic copy of
$A\,/\,\U_1\times\dots\times A\,/\,\U_n.$

Let
\begin{equation}\begin{minipage}[c]{35pc}\label{d.dim_unctb}
$J\ =\ \{i\in I\mid \dim_k(A_i)$ is uncountable$\!\}\ \subseteq\ I.$
\end{minipage}\end{equation}
If for any of $m=1,\dots,n,$ the above set
$J$ were $\!\U_m\!$-large, it is easy to see
that $A\,/\,\U_m$ would also be uncountable-dimensional.
(This does not use the countable completeness
of $\U_m;$ just the observation that if we had an $\!\aleph_1\!$-tuple
of $\!k\!$-linearly-independent elements in $A_i$ for each $i\in J,$
this would give an $\!\aleph_1\!$-tuple of elements of $A$
whose images in $A\,/\,\U_m$ would also be linearly independent.)
Then $A\,/\,\U_m$ could not be embedded in $B;$ so this does not happen.
Hence each $A\,/\,\U_m$ can be identified with a countably
complete ultraproduct of $(A_i)_{i\in I-J},$ a system of
countable-dimensional $\!k\!$-algebras.
The second paragraph of Theorem~\ref{T.nobigger}, with $\mu=\card(k)^+,$
now tells us that each $A\,/\,\U_m$ is isomorphic to one of the~$A_i.$

In the situation of statement~(a),
Lemma~\ref{L.idpt_rk} tells us that $A$ is idempotent, so by the
final parenthetical observation of Corollary~\ref{C.via_ultras+},
the $f_0$ of Theorem~\ref{T.via_A_i} is zero.
Thus, $B=f_1(A)\cong A\,/\,\U_1\times\dots\times A\,/\,\U_n,$
which we have just seen is
isomorphic to the direct product of finitely many of the $A_i.$

In the situation of~(b), $B$ will be the sum of
$f_1(A),$ as above, and $f_0(A)\subseteq Z(B).$
It is not hard to see that if an algebra $B$ is the sum of a
subalgebra $B_0\subseteq Z(B),$ and a subalgebra
$B_1$ with $Z(B_1)=\{0\},$ then $B$ is the direct product of
those two subalgebras, establishing the ``direct product''
assertion of~(b).

To get the final splitting assertion in the situation of~(a), note
that since $\U_1,\dots,\U_n$ are finitely many distinct ultrafilters,
we can partition $I$ into disjoint sets $J_1,\,\dots,\,J_n$
with $J_m\in\U_m.$
Writing $A=\prod_{J_1}A_i\times\dots\times\prod_{J_n}A_i,$
these $n$ factors have pairwise products zero, and the $\!m\!$-th
factor maps under $f$ onto the isomorphic image of $A\,/\,\U_m$ in $B.$
Thus, $f$ is, up to isomorphism, the direct product of
the $n$ canonical maps $\prod_{J_m}A_i\to\prod_{J_m}A_i\,/\,\U_m.$
As shown in Theorem~\ref{T.nobigger}, each of these
maps splits; hence so does~$f.$

The situation of~(b) is essentially the same, with the
composite map $A\to B\to B/Z(B)$ in place of~$f.$
\end{proof}

Note that if we are
given a family of idempotent algebras $A_i$ {\em not}
satisfying the hypothesis of~(a) above, there always exist
homomorphisms $f$ from $A=\prod_I A_i$ onto algebras $B$ for
which the zero-multiplication summand of statement~(b) is nonzero.
For by Lemma~\ref{L.idpt_rk}, $A$ will not be idempotent, hence $A/AA$
will be a nontrivial zero-multiplication homomorphic image of $A,$
and we can take for $B$ any countable-dimensional
homomorphic image of $A/AA.$
$(A/AA$ will itself be uncountable-dimensional,
for one can partition $(A_i)_{i\in I}$ into infinitely many
subfamilies, for each of which the finite-bound condition of~(a)
fails, and $A/AA$ maps onto the direct product of the
infinitely many zero-multiplication algebras that these yield.)

\section{Simple Lie algebras (version 2).}\label{S.simple_Lie}

What happens, in particular,
if the $A_i$ are finite-dimensional simple Lie algebras?
Will we necessarily be in situation~(a) of the above theorem?

This comes down to the question of whether, for a fixed base
field $k,$ there is a uniform bound on the idempotence ranks of
all finite-dimensional simple Lie algebras over~$k.$

At the beginning of \S\ref{S.Lie_early} we noted that a consequence of
a theorem of G.\,Brown answers that question affirmatively for $k$
algebraically closed of characteristic~$0.$
What Brown in
fact proved (in his doctoral thesis,~\cite{Brown}) is that over
{\em any} infinite field $k,$ every classical simple Lie algebra
in the sense of Steinberg~\cite{Steinberg}
has (in our language) idempotence rank~1.
(The classical simple Lie algebras in that sense comprise both
the infinite families $A_n,\,\dots\,,\,D_n$ and the exceptional
algebras, $E_6,\,\dots\,,\,G_2.)$
When $k$ has characteristic~$0,$ these are the {\em split}
simple Lie algebras in the sense of \cite{Bourbaki}, which if
$k$ is also algebraically closed are {\em all}
the finite-dimensional simple Lie algebras; so in this
case the idempotence ranks of all finite-dimensional
simple Lie algebras indeed have a common bound,~$1.$
The Bourbaki reference
that we also cited, \cite[Ch.VIII, \S13, Ex.~13(b)]{Bourbaki},
gets the same conclusion, for algebraically
closed fields of
characteristic~$0$ only, but with some additional information.

What if $k$ is, instead, the field $\R$ of real numbers?
If $L$ is a finite-dimensional simple real Lie algebra, then
$L\otimes_\R\C$ will be semisimple over $\C,$
hence a direct product of one or more simple complex Lie algebras,
hence will have idempotence rank~1 by the results quoted.
We claim that this implies that $L$
itself has idempotence rank $\leq 2.$
Indeed, every $a\in L$ can be written within
$L\otimes_\R\C$ as a bracket $[b+ic,\,d+ie]$
with $b,\,c,\,d,\,e\in L.$
Thus, $a=\mathrm{Re}(a)=[b,\,d]-[c,\,e]=[b,\,d]+[-c,\,e],$
a sum of two brackets.
(This is noted at~\cite[Corollary~A3.5, p.653]{KHH+SAM}, while
Theorem~A3.2 on the same page shows that every {\em compact} simple real
Lie algebra, i.e., every simple real Lie algebra whose Killing form is
negative definite, has idempotence rank~$1.)$
Note that the above argument
uses the fact that $\C$ has degree $2$ over $\R.$
Hence it cannot be extended to give finite bounds
on the idempotence ranks of simple Lie algebras
over most subfields $k\subseteq\C,$ e.g., $\mathbb{Q},$
since $\C$ has infinite degree over these.
(The fields for which it works, those over which an algebraically
closed field of characteristic~$0$ has finite degree, which is
necessarily~$2,$ are the real-closed fields,
the fields that algebraically ``look like''~$\R.)$

However, there is another result in the literature, less
obviously related to idempotence rank, that we can use to
get what we need in a much wider class of cases.
J.-M.\,Bois~\cite{Bois} proves, using the recently completed
classification~\cite{char>3}
of finite-dimensional simple Lie algebras $L$ over
algebraically closed fields of characteristic not~$2$ or~$3,$
that every such algebra is generated as a Lie algebra by two elements.
We shall show below, first, that such a bound on the number of
generators yields something slightly stronger than a bound
on the idempotence rank of $L,$ and, then, that for that strengthened
version of idempotence
rank, change of base field is not a problem; so that
from Bois's result on Lie algebras over algebraically closed fields,
we can get the result we need for Lie algebras over general
infinite fields.

(Notes to the reader of~\cite{Bois}: Though Theorem~A thereof
does not state the assumption that the base field is algebraically
closed, this is clear from the rest of the paper, and Bois
(personal communication) confirms that it is to be understood.
In \cite[\S1.2.2]{Bois}, the one part of that paper
where non-algebraically-closed base fields $k$ are considered,
it is shown that if $L$ is a Lie algebra
over an infinite field $k$ such that, on extending scalars
to the algebraic closure $K$ of $k,$ the resulting
Lie algebra $L\otimes_k K$ can be generated over $K$
by two elements, then $L$ can be so generated over $k.$
But we shall see from examples in \S\ref{S2.non_alg_cl} below that
in positive characteristic, a simple $L$ can yield an
$L\otimes_k K$ that is not even semisimple,
so that~\cite[\S1.2.2]{Bois} is not applicable to it; and
indeed that such an $L$ can fail to be generated by $2$ elements.
Nevertheless, the ideas of \cite[\S1.2.2]{Bois} will be used in
proving Theorem~\ref{T.x1x2} below, which states
that any such simple Lie algebra has idempotence rank $\leq 2.)$

The fact which turns statements about numbers of generators into
statements relevant to idempotence rank, part~(c) of the next lemma,
would be trivial if we were considering associative algebras.
It takes a bit more work in the Lie case, where we must
use the Jacobi identity instead of associativity, and
is false in general nonassociative algebras
(\S\ref{S2.weaker} below, last sentence).
Statements~(a) and~(b) are steps in the proof that seemed
worth recording.
These results do not require the base ring to be a field, so
we give them for general~$k.$

\begin{lemma}\label{L.gen}
Let $L$ be a Lie algebra over a commutative ring $k.$\vspace{0.2em}

\textup{(a)} \ If $U$ is a $\!k\!$-submodule of $L,$ then
$\{x\in L\mid [x,L]\subseteq U\}$ is a Lie subalgebra of
$L.$\vspace{0.2em}

\textup{(b)} \ If $V$ is a $\!k\!$-submodule of $L$ that
generates $L$ as a Lie algebra, then $[V,L]=[L,L].$\vspace{0.2em}

\textup{(c)} \ If $[L,L]=L,$ and $L$ is generated as a Lie
algebra by a set $X,$ then $L=\sum_{x\in X}\,[x,L].$
\end{lemma}

\begin{proof}
In~(a), the fact that the indicated set is
closed under the $\!k\!$-module operations
follows from the fact that $U$ is, while closure under Lie
brackets comes from the Jacobi identity:
if $[x,L]$ and $[y,L]$ are contained in $U,$ then
\begin{equation}\begin{minipage}[c]{35pc}\label{d.xyL}
$[\,[x,y],L]\ \subseteq\ [x,[y,L]\,]+[y,[x,L]\,]
\ \subseteq\ [x,L]+[y,L]\ \subseteq\ U.$
\end{minipage}\end{equation}

To get~(b), we apply~(a) with $U = [V,L].$
The Lie subalgebra described in~(a)
then contains $V,$ hence, as $V$ generates $L,$ it equals $L.$
This means that $[L,L]\subseteq [V,L];$ the opposite inclusion is clear.

To get~(c), we apply~(b) with $V$ the $\!k\!$-submodule spanned by
$X,$ and use the assumption $[L,L]=L$ to replace
$[L,L]$ in~(b) by $L.$
\end{proof}

(One can prove, more generally, a version of the above lemma
for the action of $L$ on an $\!L\!$-module $M.$
E.g.,~(c) then takes the form, ``If $L\,M=M$
and $L$ is generated as a Lie algebra by a set $X,$
then $M=\sum_{x\in X} x\,M.$'')

If we apply~(c) to the case where $L$ can be
generated by $n$ elements, the resulting conclusion,
\begin{equation}\begin{minipage}[c]{35pc}\label{d.L=[]+...}
$(\exists\ x_1,\dots,x_n\in L)
\ \ L\ =\ [x_1,L]\,+\,\dots\,+\,[x_n,L],$
\end{minipage}\end{equation}
is formally stronger than the statement that $L$ has idempotence
rank~$n:$ the idempotence rank statement allows both arguments in the
brackets giving an element $a\in L$ to vary as we vary $a,$
while~(\ref{d.L=[]+...}) fixes one argument in each bracket.
To see that it is strictly stronger, recall
that by Brown's result, many finite-dimensional
Lie algebras over fields have idempotence rank~$1;$
but no nonzero finite-dimensional Lie algebra over a field can satisfy
$(\exists\,x_1\in L)\ L=[x_1,L],$
since any $x_1$ has nontrivial centralizer, so that
$[x_1,L]$ has smaller dimension than~$L.$

Using the above lemma, Bois's Theorem~A, and a
density argument, we can now prove

\begin{theorem}\label{T.x1x2}
Let $k$ be an infinite field of characteristic not~$2$ or~$3,$
and $L$ a finite-dimensional simple Lie algebra over $k.$
Then there exist $x_1,\,x_2\in L$ such that $L=[x_1,L]+[x_2,L].$
\end{theorem}

\begin{proof}
Let $K$ be the algebraic closure of $k,$ and $L_K=L\otimes_k K.$
This will be a finite-dimensional Lie algebra over $K,$ and will inherit
from $L$ the property of being idempotent; hence as a $\!K\!$-algebra
it will have a finite-dimensional simple factor algebra $M.$
Let $q:L_K\to M$ be the canonical surjection.
Since $L$ generates $L_K$ as a $\!K\!$-algebra, $q(L)$ similarly
generates $M,$ hence, in particular, is nonzero; so, as $L$ is simple,
$q$ embeds $L$ in $M.$

Since $k$ is infinite, $L\times L$ is Zariski-dense
in $L_K\times L_K.$
(I.e., if we represent elements of the finite-dimensional
$\!K\!$-vector-space $L_K\times L_K$ in terms of coordinates
in some $\!K\!$-basis, then any polynomial function of those coordinates
which vanishes on the subset $L\times L$ vanishes everywhere.)
It follows that its image $q(L)\times q(L)\subseteq M\times M$ is
Zariski-dense in the latter space.

In what follows, let us identify $L$ with $q(L).$

By Bois~\cite[Theorem~A]{Bois}, $M$ can be generated as a
Lie algebra over $K$ by two elements.
Moreover, as noted in \cite[\S1.2.2]{Bois}, the set of generating pairs
of elements of $M$ will be a Zariski-open subset of $M\times M.$
(I.e., for every generating pair $(x_1,\,x_2)\in M\times M,$
there is a finite family of polynomials in the coordinates
of $x_1$ and $x_2$ which are nonzero at that pair, and such that
every pair at which these polynomials are nonzero is again a
generating pair.
Roughly, this is because any Lie algebra expression $f(x_1,x_2)$
has coordinates given
by polynomials in the coordinates of $x_1$ and $x_2,$ and the property
that a given list of $\dim_K(M)$ such expressions spans
$M$ over $K$ is equivalent to the nonvanishing of an appropriate
determinant in the resulting coordinate polynomials.
Cf.~\cite[\S1]{x1Adots}.)
The nonempty Zariski-open set of generating pairs must meet the
Zariski-dense set $L\times L,$ which means
that there exist $x_1,\,x_2\in L$ which
generate $M$ as a Lie algebra over $K.$
Hence by Lemma~\ref{L.gen}(c),
\begin{equation}\begin{minipage}[c]{35pc}\label{d.M=}
$M\ =\ [x_1,\,M]\,+\,[x_2,\,M].$
\end{minipage}\end{equation}

We claim that this implies
\begin{equation}\begin{minipage}[c]{35pc}\label{d.L=}
$L\ =\ [x_1,\,L]\,+\,[x_2,\,L].$
\end{minipage}\end{equation}

To show this, let $B$ be a basis of $K$ as a $\!k\!$-vector-space.
Then the Lie algebra $L_K=L\otimes_k K,$ under the adjoint
action of its sub-$\!k\!$-algebra $L,$ is the direct sum of the
sub-$\!L\!$-modules $L\otimes b$ $(b\in B),$
each of which is isomorphic to $L$ as an $\!L\!$-module,
via the map $x\mapsto x\otimes b.$
Since $L$ is simple as a Lie
algebra, it is a simple module over itself under the
adjoint action, hence since $L_K$ is a direct sum of copies
of that simple $\!L\!$-module, so is its homomorphic image $M.$
As $M$ is a direct sum of simple $\!L\!$-modules, $L$ is a direct
summand therein.
Applying to~(\ref{d.M=}) an $\!L\!$-module projection
of $M$ onto $L,$ we get~(\ref{d.L=}), as required.
\end{proof}

For some results on particular elements $x_1$ and $x_2$
in split simple Lie algebras $L$ such that $L=[x_1,L]+[x_2,L],$
and related matters, see \cite{NN_split}.

Theorem~\ref{T.simple}(a) and Theorem~\ref{T.x1x2} together
give the desired result on infinite products of simple Lie algebras
in characteristics $\neq 2,\,3.$
We also record the weaker statement
that follows from Theorem~\ref{T.simple}(b)
for characteristics $2$ and $3$
(where there is as yet insufficient structure theory to
say whether a result like that of~\cite{Bois} holds).

\begin{theorem}\label{T.simple_Lie}
Suppose $k$ is an infinite field, and $f: \prod_I A_i\to B$ a
surjective homomorphism from a direct product of
finite-dimensional simple Lie algebras
to a finite-dimensional Lie algebra~$B.$
Then

\textup{(a)} \ If $\mathrm{char}(k)\neq 2$ or~$3,$
$B$ is isomorphic to the direct product of finitely many
of the $A_i,$ and the map $f:A\to B$ splits.

\textup{(b)} \ If $\mathrm{char}(k)= 2$ or $3,$ one can at
least say that $B$ is isomorphic to a direct product of
finitely many of the $A_i$ and an abelian Lie algebra.
In that case, the composite map $A\to B\to B/Z(B)$ splits.\qed
\end{theorem}

A few notes on the general concept of idempotence rank:
By Theorem~8.4.5 of \cite{PMC_further},
every finite-dimensional {\em simple associative} algebra has a unit,
and so has idempotence rank~$1.$
On the other hand, we will give in~\S\ref{S2.idrk}
examples of arbitrarily large finite
idempotence rank in {\em non-simple} finite-dimensional
Lie and associative algebras, and in simple finite-dimensional
{\em non}\/-Lie {\em non}\/associative algebras;
while in \S\ref{S2.idrk_neq}, we will give
an example of a finite-dimensional (non-Lie) algebra
whose idempotence rank changes under change of base field.
(This is the phenomenon, the possibility of which
prevented us from using Brown's
result to get Theorem~\ref{T.simple_Lie} over a
general field of characteristic~$0.)$

\section{Continuity in the product topology.}\label{S.topology}

Any infinite product $A=\prod_I A_i$ of sets has a natural topology,
the product of the discrete topologies on the $A_i.$
If the $A_i$ have group structures, and $f:A\to B$ is a
homomorphism into another group $B,$ it is not hard to show
that $f$ is continuous in the product topology on $A$ and the
discrete topology on $B$ if and only if it factors through the
projection of $A$ onto a finite
subproduct $A_{i_1}\times\dots\times A_{i_n}.$

Indeed, ``if'' is immediate.
To see ``only if'', note that by the discreteness of $B,$
$\ker(f)$ is open.
As an open set containing the identity element, it must
contain the intersection of the inverse images of neighborhoods
of the identity elements $e_i$ of finitely many of the $A_i.$
So a fortiori, it contains the intersection of the inverse
images of the trivial subgroups of those $A_i;$ which is the
kernel of the projection to their product, so $f$ factors through
that projection.

Let us briefly note when this continuity condition holds
in the results of the preceding sections.

It is not hard to see that it can never hold if $f$ involves
a factorization through a {\em nonprincipal} ultrafilter.

When $f$ is the sum of a map $f_1$ that factors
through finitely many of the $A_i$
(corresponding to finitely many {\em principal} ultrafilters), and
a possibly nonzero perturbing map $f_0$ into $Z(B),$ then the
continuity of $f$ depends on the continuity of $f_0,$ which in
general cannot be expected: as noted in \S\ref{S.intro},
such maps tend to be ``unruly''.
However, the effect of $f_0$ disappears if we compose
$f$ with the factor map $\pi: B\to B/Z(B).$
Summarizing the consequences of these considerations, we have

\begin{proposition}\label{P.continuity}
In the preceding results of this note, the map $f$ will
be continuous in the product topology on $A$ and the
discrete topology on $B$ \textup{(}equivalently, will
factor through the projection to a finite subproduct
of the $A_i)$ in the situations of the following results whenever
$\card(I)$ is less than every
measurable cardinal $>\nolinebreak\card(k):$
Theorem~\ref{T.simple<} \textup{(}where that restriction on $I$
is already assumed\textup{)}, Theorem~\ref{T.simple}\textup{(a)}, and
Theorem~\ref{T.simple_Lie}\textup{(a)}.

Under the same assumption on $\card(I),$
the composite map $\pi\,f:A\to B\to B/Z(B)$ will be continuous in the
situations of Theorem~\ref{T.via_A_i}, Theorem~\ref{T.solv2},
Theorem~\ref{T.nilp}, Theorem~\ref{T.simple}\textup{(b)}, and
Theorem~\ref{T.simple_Lie}\textup{(b)}.
\end{proposition}

\section{Some tangential notes.}\label{S.tangents}

We record here further observations on the
material introduced in the preceding pages, which were not
needed for the results developed there.
(The subsections of this section are independent of one another.)

\subsection{Almost direct factors, and Boolean algebras}\label{S2.Boole}
Recall that for an {\em associative unital} algebra $A,$
an almost direct decomposition is the same as a
direct product decomposition (because $Z(A)=\{0\}),$ and that in this
situation, such decompositions are in
bijective correspondence with the central idempotent elements of $A.$
The set of such central idempotents, and hence the partially
ordered set of almost direct factors of $A,$ forms a Boolean algebra.
Will the same be true of the partially ordered set
of almost direct factors in a general algebra $A$?

Below, we obtain a positive answer when $A$ is idempotent
(or satisfies a slight weakening of that condition),
then a counterexample in the absence of that assumption.

\begin{proposition}\label{P.boole}
Let $A$ be an idempotent algebra, or more generally, an algebra
satisfying
\begin{equation}\begin{minipage}[c]{35pc}\label{d.AA+Z}
$A\,=\,AA+Z(A).$
\end{minipage}\end{equation}

Then the almost direct factors of $A$
form a Boolean algebra, with zero element $Z(A),$ unit element $A,$
the join of $B$ and $C$ given by the sum $B+C,$ and
the meet given by the intersection $B\cap C,$ which is
also equal to $BC+Z(A)$ and to $CB+Z(A).$
\end{proposition}

\begin{proof}
Let us write $a\mapsto\overline{a}$ for the quotient map $A\to A/Z(A).$
Any almost direct decomposition $A=B+B'$ is determined by the induced
direct product decomposition
$\overline{A}=\overline{B}\times\overline{B'}$
(cf.\ remarks following Definition~\ref{D.almost}), hence by the
projection operator of $\overline{A}$ onto $\overline{B}.$
We shall prove below that under the present
hypotheses, the projection operators so induced by any two
almost direct decompositions $A=B+B'$ and $A=C+C'$ commute, and
that the image $\overline{B}\cap\overline{C}$ of their
composite corresponds to an almost direct factor
$B\cap C=BC+Z(A)=CB+Z(A)$ of $A,$ with complement $B'+C'.$
Now a set of pairwise commuting projection operators
(i.e., idempotent endomorphisms) on any
abelian group generates a Boolean algebra of such operators,
with the meet and join of operators $e$ and $f$ (given by
$ef$ and $e+f-ef$ respectively)
corresponding to the intersection and the sum
of the image subgroups; so these results will prove our claims.

Given almost direct decompositions $A=B+B'$ and
$A=C+C',$ let us multiply these two equations together and add $Z(A).$
By~(\ref{d.AA+Z}),
this yields $A=BC+BC'+B'C+B'C'+Z(A),$ which we can rewrite
\begin{equation}\begin{minipage}[c]{35pc}\label{d.A=BC+}
$A\,=\,(BC+Z(A))\,+\,(BC'+Z(A))\,+\,(B'C+Z(A))\,+\,(B'C'+Z(A)).$
\end{minipage}\end{equation}
(Since $A$ is not assumed associative or Lie, we do not yet know
that the summands in~(\ref{d.A=BC+}) are ideals of $A,$ only that
they are $\!k\!$-submodules.)

Let us verify first that the decomposition
$a=a_{BC}+a_{BC'}+a_{B'C}+a_{B'C'}$ of an element $a\in A$
arising from~(\ref{d.A=BC+}) is unique modulo $Z(A).$
For this, it will suffice to show that in any decomposition of $0,$
\begin{equation}\begin{minipage}[c]{35pc}\label{d.dec_0}
$0\,=\,z_{BC}+z_{BC'}+z_{B'C}+z_{B'C'}$
\end{minipage}\end{equation}
into summands in the above four $\!k\!$-submodules,
all of these summands must lie in $Z(A).$
Now since $B$ and $C$ are ideals, the term $BC$
in the first summand of~(\ref{d.A=BC+}) is contained
in both $B$ and $C,$ hence that summand $BC+Z(A)$
annihilates $B'$ and $C',$ hence
annihilates all the summands in~(\ref{d.A=BC+}) other than itself;
so if the summand $z_{BC}$ of~(\ref{d.dec_0}) does not lie in
$Z(A),$ i.e., does not annihilate all of $A,$ this can only be
because it fails to annihilate the first summand of~(\ref{d.A=BC+}).
But all the other summands on the right-hand side of~(\ref{d.dec_0}) do
annihilate the first summand of~(\ref{d.A=BC+}), as does
the left-hand term, $0;$ so $z_{BC}$ must also.
This completes the verification that it lies in $Z(A);$
and by the same argument, so do all the terms of~(\ref{d.dec_0}),
as claimed.
Hence, passing to quotients modulo $Z(A),$ the decomposition
\begin{equation}\begin{minipage}[c]{35pc}\label{d.barA=}
$\overline{A}\,=\,\overline{BC}+\overline{BC'}+
\overline{B'C}+\overline{B'C'}.$
\end{minipage}\end{equation}
is a direct product decomposition of $\!k\!$-modules.

Using again the same kind of reasoning,
note that when we decompose an element of $A$ by~(\ref{d.A=BC+}),
the component in each summand whose expression in~(\ref{d.A=BC+})
involves $B$ annihilates $B',$ and inversely.
Hence given such a decomposition $a=a_{BC}+a_{BC'}+a_{B'C}+a_{B'C'},$
the expression $a=(a_{BC}+a_{BC'})+(a_{B'C}+a_{B'C'})$ decomposes $a$
into an element annihilating $B',$ i.e., a member of $B,$
and an element annihilating $B,$ i.e., a member of $B'.$
But the decomposition of $a$ coming from the relation $A=B+B'$ is
unique up to summands in $B\cap B'=Z(A);$
hence the idempotent endomorphism
of $\overline{A}$ given by projection on the first summand
in $\overline{A}=\overline{B}+\overline{B'}$ must coincide with the
projection of~(\ref{d.barA=}) onto the sum of its first and
second summands.
Similarly, the idempotent endomorphism of $\overline{A}$ arising
from the decomposition $A=C+C'$ must be the projection
of~(\ref{d.barA=}) onto the sum of its first and third summands.
These two projections commute, since their product in either
order is the projection of~(\ref{d.barA=}) onto its first summand.

Since the range of the product of two commuting idempotent
endomorphisms of an abelian group is
the intersection of their ranges, we have
$\overline{BC}=\overline{B}\cap\overline{C}.$
Taking inverse images in $A,$ this gives
\begin{equation}\begin{minipage}[c]{35pc}\label{d.BC+=cap}
$BC+Z(A)\,=\,B\cap C,$
\end{minipage}\end{equation}
as claimed;
and by symmetry we likewise have $CB+Z(A)=B\cap C.$

The equality~(\ref{d.BC+=cap}) shows that $BC+Z(A)$ is an ideal;
we must still verify that it is an almost direct factor in $A.$
We claim that it and $(BC'+Z(A))+(B'C+Z(A))+(B'C'+Z(A))$
are each other's two-sided annihilators.
We have seen that they annihilate each other.
On the other hand, by the method of reasoning used immediately
after~(\ref{d.dec_0}), if an element $a$
annihilates $BC+Z(A),$ the first component of
a decomposition of $a$ as in~(\ref{d.A=BC+}) lies
in $Z(A);$ and since $Z(A)$ also lies in the other three summands
of~(\ref{d.A=BC+}), $a$ will lie in the sum of those three summands.
So the two-sided annihilator of $BC+Z(A)$ is indeed
$(BC'+Z(A))+(B'C+Z(A))+(B'C'+Z(A)).$
That $BC+Z(A)$ is likewise the two-sided annihilator of that sum is
shown in the same way.
This completes our proof.
\end{proof}

The above result covers not only the case where $A$ is idempotent,
but the opposite extreme, where $A$ has zero multiplication,
since then $A=Z(A)=AA+Z(A).$
(In that case, our Boolean algebra is trivial.)
But let us now show that when $A\neq AA+Z(A),$
the conclusion of Proposition~\ref{P.boole} need not hold.

Let $A=\R^2\times \R,$ with
multiplication $(v,a)*(w,b)=(0,v\cdot w),$ where $v\cdot w$
is the dot product of vectors in $\R^2.$
Then $Z(A)=\{0\}\times\R,$ and
for every one-dimensional subspace $V\subseteq\R^2,$
we have the almost direct decomposition
$A=(V\times\R)+(V^\perp\times\R),$ where
$(~)^\perp$ denotes orthogonal complement in $\R^2.$
Thus, the almost direct factors lying strictly between
$Z(A)$ and $A$ form an infinite set of pairwise incomparable
elements $V\times\R$ (though for each such element,
only one of the others is its ``complementary almost direct factor''
as we have defined the term).
Hence the partially ordered set of almost
direct factors of $A$ is not a Boolean algebra.

The above algebra $A$ is, incidentally, associative, since
all three-fold products are zero.

\subsection{Weakening the definition of an almost direct decomposition}\label{S2.weaker}
In our definition of an almost direct decomposition $A=B+B'$ of an
algebra $A,$ the condition that $B$ and $B'$ be ideals can
be formally weakened to say that they are subalgebras.
For the latter condition on $B$ says that it is closed
under left and right multiplication by $B,$ and since it annihilates
$B',$ it is trivially closed under left and right multiplication by
that subalgebra; hence it closed under left and right multiplication by
$B+B'=A.$

If we ask whether it is enough to assume that $B$ and $B'$
are $\!k\!$-submodules summing to $A,$ each of
which is the other's two-sided annihilator, the answer is mixed.
If $A$ is associative, we can still conclude that they
will be almost direct factors.
For since $B$ annihilates $B'$ on both sides, associativity implies
that $BB$ does the same; hence it is contained in the annihilator
of $B',$ namely $B,$ proving that $B$ is a subalgebra.
We get the same conclusion if $A$ is a Lie algebra: the
Jacobi identity shows that $[B',[B,B]]\subseteq [B,[B',B]]=\{0\},$
hence $[B,B]$ is contained in the annihilator $B$ of $B'.$

But for general $A,$ the corresponding statement is false.
Indeed, for any $k,$ let $A$ be the $\!k\!$-algebra
which is free as a $\!k\!$-module on
two elements $x$ and $y,$ with multiplication given by
\begin{equation}\begin{minipage}[c]{35pc}\label{d.xy=0}
$xy=yx=0,\qquad xx=y,\qquad yy=x.$
\end{minipage}\end{equation}
Then clearly, $kx$ and $ky$ are each other's two-sided annihilators,
and sum to $A,$ but are not subalgebras.

Even if one of $B,\ B'$ is a subalgebra,
the other may not be, as we can see by replacing the relation
$xx=y$ in~(\ref{d.xy=0}) by $xx=x,$ while leaving the other
relations unchanged.

Incidentally, taking $X=\{x\}$ in the algebra defined
by~(\ref{d.xy=0}), we find that the analog of Lemma~\ref{L.gen}(c)
fails $(X$ generates $A,$ but $XA\neq A),$
showing that that result does not
hold for general $\!k\!$-algebras.

\subsection{``Early'' ultrafilters}\label{S2.early}
Just as many calculus texts come in two versions, ``early
transcendentals'' and ``late transcendentals'', so the
development of \S\S\ref{S.prelim}-\ref{S.prod_gen}
has an alternative version, in
which we obtain our ultrafilters early, before the
``either/or'' conditions such as Lemma~\ref{L.A1A2->B}(i)
by which we summoned them in our present development.

In such a development, one would associate to any map $f$ from a product
of nonempty sets $A=\prod_I A_i$ to a set $B$
the family $\F_f$ of $J\subseteq I$ such that $f$ factors through
$\prod_J A_i.$
This turns out to be a filter, the largest filter $\F$
such that $f$ factors through $A/\F.$
Any filter is an intersection of ultrafilters; let us call the set of
ultrafilters containing $\F_f$ the ``support'' of $f.$
One verifies that $f$ factors through
the natural map $A\to\prod_{\U\supseteq\F_f} A/\U.$
Finally, bringing in the assumption that the $A_i$ and $B$ are algebras
and $f$ a surjective homomorphism, one can use the argument of
Proposition~\ref{P.via_ultra} to show that if $B$ satisfies
the conditions of Lemma~\ref{L.A1A2->B}, then
the support of $f$ is a singleton $\{\U\},$
while the argument of Proposition~\ref{P.via_ultras/Z} shows that
if $B$ satisfies the weaker property of chain condition on
almost direct factors, then
$\pi f$ has support in a finite set of ultrafilters.

The proofs of the propositions mentioned used the fact that in a direct
product of algebras,
elements with disjoint support have trivial product.
One might get similar results on direct products of groups
(or even monoids) using the fact that in a direct product
of these, elements with disjoint supports commute.
(This is suggested, of course, by the way the brackets
of Lie algebras arise from the noncommutativity of Lie groups.)
We leave this for the interested reader to investigate.

A very different ``early ultrafilters'' approach is
taken in \cite[\S3]{prod_Lie2}.

\subsection{On idempotence rank, and related functions.}\label{S2.r(M)}
In examining the properties of the idempotence rank function
on idempotent algebras, it is helpful to look at
a more general version of that situation.
For simplicity, let $k$ be a field.
Consider any $\!4\!$-tuple
$(A,B,C,m),$ where $A,\ B$ and $C$ are $\!k\!$-vector-spaces,
and $m$ is a surjective linear map $A\otimes_k B\to C.$
(Thus, $m$ gives the same information as a $\!k\!$-bilinear
map $A\times B\to C$ whose image spans $C.)$
Let us define
\begin{equation}\begin{minipage}[c]{35pc}\label{d.maxrank}
$\mbox{max-rank}(m)\,=\,
\sup_{c\in C}\,\inf_{t\in m^{-1}(c)}\,\mathrm{rank}(t),$
\end{minipage}\end{equation}
where $\mathrm{rank}(t)$ denotes the rank of $t$ as a member
of the tensor product $A\otimes_k B,$ i.e., the minimum number
of decomposable tensors $a\otimes b$ that must be summed to get $t.$
We see that when $A=B=C$ is the underlying vector space of an
idempotent $\!k\!$-algebra $A,$ and $m$ the map corresponding to the
multiplication of $A,$
then $\mbox{max-rank}(m)$ is in fact $\idrk(A).$

The function $\mbox{max-rank}(m)$ has a family resemblance to
the function $r(M)$ introduced in \cite{r(M)}
for a subspace $M$ of a tensor product $A\otimes_k B,$ and defined by
\begin{equation}\begin{minipage}[c]{35pc}\label{d.r(M)}
$r(M)\,=\,\inf_{t\in M-\{0\}}\,\mathrm{rank}(t).$
\end{minipage}\end{equation}
The contexts of the two definitions are essentially the same:
what we are given in each is equivalent to a short exact sequence
$0\to M\to A\otimes_k B\to C\to 0$ of $\!k\!$-vector-spaces,
with middle term a tensor product.
However, neither of these invariants of that short exact sequence
seems to be expressible in terms of the other.
In fact, if our vector spaces are finite-dimensional, we can form the
dual short exact sequence $0\to C^*\to A^*\otimes_k B^*\to M^*\to 0,$
and look at the same two invariants for it, getting, altogether,
four invariants from our original sequence, none of which seems
to be expressible in terms of the others.

As noted in \cite{r(M)}, $r(M)$ can decrease,
but not increase, under extension of base field;
for when we make such an extension, the set of elements
over which the infimum of~(\ref{d.r(M)}) is taken is enlarged,
while the rank-function on elements lying in the original
tensor product remains unchanged.
(If we take for $M$ the kernel of the
map $A\otimes_k A\to A$ corresponding to the multiplication
operation of an algebra, then a decrease in $r(M)$
from a value $>1$ to $1$ under base extension from $k$ to $K$
means that from a $\!k\!$-algebra $A$ without zero divisors, we get
a $\!K\!$-algebra $A\otimes_k K$ with zero divisors.
The reverse cannot happen, of course.)

Since the definition~(\ref{d.maxrank})
of $\mbox{max-rank}(m)$ involves both a supremum and
an infimum, that function can potentially
increase or decrease under base extension.
The possibility of its decreasing is what made it impossible
for us to go from Brown's result showing that $\idrk(L)=1$ for $L$ a
finite-dimensional simple Lie algebra over $\C$
to the corresponding statement for subfields of $\C.$
In \S\ref{S2.idrk_neq} we will see examples
of finite-dimensional, idempotent (but non-Lie, nonassociative,
non-simple) algebras $A$ whose idempotence ranks do increase
and decrease under base extensions.
We do not know whether either can happen
when $A$ is a simple Lie algebra.

What we used in \S\ref{S.simple}
(in conjunction with the results of~\cite{Bois}),
instead of the unsuccessful approach indicated above, was
an argument via what might
be called the ``one-variable-constant idempotence rank function'',
the least $n$ such that~(\ref{d.L=[]+...}) holds.
In the case of general $\!k\!$-algebras, where left and right
multiplication are not equivalent, we could call the version with
the constant factors on,
say, the left the ``left-constant idempotence rank'':
\begin{equation}\begin{minipage}[c]{36pc}\label{d.lconst}
$l\!$-const-$\!\idrk(A)\,=\,\inf\,
\{n\mid (\exists\,x_1,\dots,x_n\in A)\ A=x_1A+\dots+x_nA\}
\,=\,\inf\,\{\dim_k(V)\mid A=VA\}.$
\end{minipage}\end{equation}

This function is examined in~\cite{x1Adots}.

\subsection{Other literature on homomorphisms from infinite products.}\label{S2.other}
Restrictions on homomorphisms from infinite direct products
to ``small'' objects have been noted in other areas of algebra.

In \cite[Corollary~9]{BHN+SY}, it is shown that
a homomorphic image of a direct power of
a finite nonabelian simple group $G,$ if countable, must be finite;
general finite groups $G$ with that property are investigated in
\cite{MB} and \cite{MJB}.
This situation has a similar
flavor to that of the present note; e.g., note that
simple nonabelian groups satisfy the analogs of
the conditions of our Lemma~\ref{L.A1A2->B}.
(We remark, however, that the groups characterized in the above
papers all have trivial centers.
Perhaps if one considers groups
with nontrivial centers, analogs of the results of this note showing
that homomorphisms $\prod A_i\to B$ acquire stronger properties on
composition with the natural map $B\to B/Z(B)$ will turn up.)

An area of investigation with a different flavor begins with the result
of~\cite{Specker}, that every homomorphism of abelian
groups $\mathbb{Z}^\mathbb{N}\to\mathbb{Z}$ factors
through the projection onto finitely many coordinates.
It can be deduced from this that the same factorization property holds
for homomorphisms from any countable product of abelian groups
$\prod_{i\in\mathbb{N}}A_i$ to $\mathbb{Z};$ this is
expressed by saying that $\mathbb{Z}$ is a {\em slender} abelian group.
More generally, slenderness has been studied in abelian monoids,
in modules over general rings, and in
objects of general preadditive categories.
Note that for these abelian groups,
abelian monoids, etc., unlike the algebras of this note
and unlike nonabelian groups, any finite family of
morphisms can be added; hence in mapping
a finite product $A_1\times\dots\times A_n$ to an object $B,$
one can form sums of homomorphisms $A_i\to B.$
Thus, the restrictions that turn out to
hold on homomorphisms from infinite products of these objects
cannot arise from restrictions on homomorphisms from finite products,
like those of our Lemma~\ref{L.A1A2->B}, but must come in
in a more mysterious way; roughly, it seems, from completeness-like
properties of infinite products, which cannot be duplicated
in a ``slender'' $B.$
Examples of abelian groups that are {\em not} slender include
all abelian groups with torsion, all nonzero injective abelian
groups, the additive group of $\!p\!$-adic integers, and, of course,
$\mathbb{Z}^\mathbb{N}.$
For a sampling of work in this area, see,
for abelian groups,~\cite[\S94]{FuchsII} and~\cite{SS+OK},
for modules,~\cite{DA} and~\cite{B+K1}, and for
abelian monoids and objects
of preadditive categories,~\cite{RD:mnd} and~\cite{RD:bk}.

Related conditions have been considered on nonabelian
groups, in some cases again defined in terms of homomorphisms
from direct products~\cite{RG}~\cite{F44}, in others,
in terms of homomorphisms from certain completions of free
products \cite{Higman} \cite{KEda}~\cite{ncmSpecker}.

We remark that for abelian groups and other structures whose
morphisms can be added, the class of slender objects
would not change if, in the definition, we restricted attention
to {\em surjective} homomorphisms $\prod A_i\to B,$
since if $f:\prod A_i\to B$ is a nonsurjective
map witnessing the failure of $B$ to be slender, there
is an obvious surjective map $B\times\nolinebreak\prod A_i\to B$
which does the same.
A similar observation applies to the version of slenderness in
nonabelian groups defined using the ``complete free products''
of~\cite{KEda} --
but not to the one defined using direct products~\cite{RG}.
Hence if one defines a condition like slenderness for nonabelian
groups, but based on {\em surjective} maps from direct products, one can
expect to find a larger class of examples than the ordinary slender
groups, and probably techniques and results close to those of this note;
cf.\ next to last paragraph of \S\ref{S2.early}.
(It is not clear to us whether the class of groups defined
similarly in terms of maps from the ``unrestricted free
products'' of \cite{Higman},~\cite{ncmSpecker}
would similarly grow if one imposed this condition
only on {\em surjective} maps from those groups.)

In~\cite{E+F+M}, some implications among conditions
on homomorphisms $A^I\to B$ are studied for algebras $A$ and $B$
in the general sense of
universal algebra, the cases of slender abelian groups on the
one hand, and of discriminator algebras on the other, being noted.

In this section we have, for simplicity, limited the results quoted
to the countable-index-set case; though in the works cited,
what is in question is generally whether the
index set is smaller than all uncountable measurable cardinals.

\section{Examples.}\label{S.examples}

In earlier sections, we noted some examples in passing.
Here we give further, mostly lengthier examples, for which we did not
want to interrupt the development of our earlier results.

As in \S\ref{S.tangents}, the subsections below are independent
of one another.
The only dependence on that section is
that \S\ref{S2.idrk_neq} below assumes~\S\ref{S2.r(M)} above.

\subsection{Idempotent algebras with $Z(A)\neq\{0\}\!$}\label{S2.ZinIdpt}
We noted following Proposition~\ref{P.k->K} that
a {\em unital} algebra $A$ necessarily satisfies $Z(A)=\{0\}.$
Is the same true of {\em idempotent} algebras -- perhaps
subject to some additional conditions?

An easy example shows that this need not even hold in
finite-dimensional algebras of idempotence rank~$1.$
Let $\mathbb{H}$ be the $\!\R\!$-algebra of
quaternions, let $\mathrm{Im}: \mathbb{H}\to\mathbb{H}$ be the
``imaginary part'' map,
$a+bi+cj+dk\mapsto bi+cj+dk,$ and let $A$ be $\mathbb{H}$ under the
nonassociative multiplication $x*y=\mathrm{Im}(x)\,\mathrm{Im}(y)$
(where the right-hand side is evaluated
using the ordinary multiplication of $\mathbb{H}).$
Note that if we call the real and imaginary parts of an
element of $\mathbb{H}$ its ``scalar'' and ``vector'' components,
then $x*y$ has for scalar component the negative of the dot product of
the vector components of $x$ and $y,$
and for vector component the cross product of those same vectors.
Now it is not hard to see geometrically that for any scalar $a$ and
vector $bi+cj+dk,$ one can find two vectors with dot product $-a$ and
cross product $bi+cj+dk.$
This gives the asserted idempotence of our algebra.
On the other hand, clearly, $Z(A)=\R\neq\{0\}.$

One can get an infinite-dimensional
example, again of idempotence rank~$1,$ that is
associative and commutative:
Let $V$ be any commutative valuation ring with
nondiscrete valuation; thus, its maximal ideal $\mathbf{m}$
is idempotent of idempotence rank $1.$
Take a nonzero element $x\in\mathbf{m},$
and let $A=\mathbf{m}/x\,\mathbf{m}.$
Then $A$ has idempotence rank $1,$ but the image of $x$ lies in $Z(A).$

We shall also see in \S\ref{S2.idrk} below, where we give
examples of finite-dimensional idempotent associative algebras
$A$ of arbitrarily large finite idempotence rank,
that such algebras can have $Z(A)\neq\{0\}.$

\subsection{On the chain condition on almost direct factors.}\label{S2.CC}
Our next example will show that a finitely generated (unital
or nonunital) associative
algebra need not satisfy the chain condition on almost direct factors;
and thus that the Boolean algebra of central idempotents of
a finitely generated unital associative algebra can be infinite.

We begin by constructing a family of unital associative algebras $A_i$
$(i=0,1,\dots)$ over any field $k.$
For each $i,$ let $A_i$ be presented by three
generators $x,\,y,\,z,$ and the infinite family of relations
\begin{equation}\begin{minipage}[c]{35pc}\label{d.xynz}
$xy^nz\,=\,\left\{ \begin{array}{ll}
1\quad\mbox{if}\ \ n=i\\
0\quad\mbox{otherwise}
\end{array}
\right.$
$(n=0,1,\dots).$
\end{minipage}\end{equation}

If we regard~(\ref{d.xynz}) as a system of ``reduction rules''
for expressions in $x,\,y,\,z,$ we find that, in the terminology
of~\cite{<>}, these have no ``ambiguities'' (roughly, there is
no way to write down a word having
subwords $xy^mz$ and $xy^nz$ which ``overlap'',
and so force us to worry whether the two competing reductions may fail
to lead ultimately to the same expression).
Moreover, application of one of these reduction rules to a
word in $x,$ $y$ and $z$ yields at most a shorter word, so
the process of recursively applying these rules always terminates.
Hence, by \cite[Theorem~1.2]{<>}, each $A_i$ has for $\!k\!$-basis
the set of words in $x,\,y,\,z$ (including the empty word, $1,$
since at the moment we are considering unital $\!k\!$-algebras)
having no subwords of the form $xy^nz.$
In particular, the empty word $1$ belongs to this
basis, so $1\neq 0$ in each of these algebras.

Now within the product algebra $\prod_{i=0,1,\dots} A_i,$
let $x$ be the element whose coordinate in each
$A_i$ is the element $x\in A_i,$ define $y$ and $z$ analogously,
and let $A\subseteq \prod A_i$ be the
nonunital subalgebra generated by these three elements.
Then for each $n,$ the element $xy^nz\in A$
will have $1$ in the $\!n\!$-th
coordinate and $0$ in all others, and so be a central idempotent.
It follows that the elements
$xz,\ xz+xyz,\,\dots\,,\ \sum_{m=0}^n xy^mz,\,\dots$ constitute
an infinite ascending chain of central idempotents,
yielding an infinite ascending chain of almost direct factors.

If, instead, we take for $A$ the {\em unital} algebra generated
by these same three elements, we get
a finitely generated unital associative algebra whose
Boolean algebra of central idempotents is infinite.

We remark that in~(\ref{d.xynz}), for conceptual simplicity, we used
algebras $A_i$ such that the sets $\{n\mid xy^nz=1\}$ were singletons;
but for every subset $J$ of the natural numbers, the same structure
result applies to the algebra $A_J$ obtained by setting $xy^nz$ to
equal $1$ if $n\in J,$ and $0$ otherwise.
Applying the above construction to appropriate families of
these algebras, we can get a $\!3\!$-generator algebra $A$
whose Boolean algebra of central idempotents is any
countable Boolean algebra.

\subsection{The need for $k$ to be a field in Lemma~\ref{L.AAcapC}.}\label{S2.nonfd}

In Lemma~\ref{L.AAcapC}, assuming $k$ a field, we gave a necessary
and sufficient condition for a homomorphism $f: A\to B$ of
$\!k\!$-algebras to be approximable modulo $Z(B)$ by a homomorphism
$f_1: A\to B$ annihilating a given ideal $C\subseteq A.$
The following example shows that the condition
given there fails to be sufficient if, instead, $k$ is any
integral domain that is not a field (or more generally,
if $k$ is a commutative ring that is not von Neumann regular).

\begin{lemma}\label{L.nonfield}
Let $k$ be a commutative ring having an element $c$ such
that $c\notin c^2k.$

Let $B$ be the free $\!k/c^2k\!$-module
on one generator $x,$ with the \textup{(}associative,
commutative\textup{)} multiplication defined by $x^2=cx;$
let $A_0$ be the free $\!k/c^2k\!$-module
on one generator $y,$ with the zero multiplication;
let $A=B\times A_0$ as $\!k\!$-algebras,
and let $f:A\to B$ be the projection onto the first factor.
Then the $\!k\!$-submodule $C$ of $A$ generated by $cx-cy$ is an
ideal that satisfies condition~\textup{(ii)} of Lemma~\ref{L.AAcapC},
but not condition~\textup{(i)}.
\end{lemma}

\begin{proof}
Observe first that $cx-cy\in Z(A),$ hence the $\!k\!$-submodule
$C$ that it generates is indeed an ideal, and its image under $f$
lies in $Z(B).$

Note also that if an element $d(cx-cy)$ $(d\in k)$ of this ideal lies
in $AA,$ then its $\!A_0\!$-component $dcy$ must lie
in $A_0A_0=\{0\},$ hence since the subrings $A_0=ky$ and $B=kx$
are isomorphic as $\!k\!$-modules, we also have $dcx=0,$
so the given element $d(cx-cy)$ is zero.
Thus $C\cap AA=\{0\},$ so our example satisfies
condition~\textup{(ii)} of Lemma~\ref{L.AAcapC}.

Now suppose there were a homomorphism
$f_1$ as in condition~(i) of that lemma.
Lemma~\ref{L.fgh}(ii) tells us that $f-f_1$ annihilates $AA,$
which contains $x^2=cx,$ so $f_1$ agrees with $f$ at $cx,$
i.e., it fixes that element.
But by assumption, $f_1$ annihilates $cx-cy\in C,$ so
we must also have $f_1(cy)=cx.$

Now writing $f_1(y)=ax\in B,$ the above equation becomes $c(ax)=cx.$
Applying this twice, we get
$cx=acx=a^2cx=a^2x^2=(ax)^2=f_1(y)^2=f_1(y^2)=f_1(0)=0,$
a contradiction.
(Intuitively, the algebra structures on $kx$ and $ky$ are too different
for there to be a nice choice of $f_1$ annihilating $cx-cy.)$
\end{proof}

\subsection{Unbounded idempotence rank.}\label{S2.idrk}
Lemma~\ref{L.idpt_rk} tells us that a product $A=\prod_I A_i$ of
finite-dimensional idempotent algebras will fail to be idempotent
if the idempotence ranks of those algebras are unbounded; but
we have seen that in most characteristics, finite-dimensional
{\em simple Lie} algebras all have idempotence rank $\leq 2,$ and
even in the
remaining two characteristics, one may hope that the the same is true.
Can we get any examples of finite-dimensional idempotent algebras with
arbitrarily large finite idempotence ranks?

We give below three classes of such examples:
for associative algebras, for Lie algebras, and for
(nonassociative non-Lie) simple algebras, respectively.
(Our descriptions of the first two constructions
also record the fact that $Z(A)$ is nontrivial, giving
examples mentioned in the
last paragraph of \S\ref{S2.ZinIdpt}.

\begin{lemma}\label{L.ae11b}
For any field $k$ and positive integer $i,$ let $A$ be
the $\!k\!$-algebra with underlying vector space
the space $M_i(k)$ of $i\times i$ matrices over $k,$ and
multiplication ``$*\!$'' expressed, in terms of the ordinary
multiplication of these matrices, by $a*b=a\,e_{11}\,b.$

Then $A$ is an associative idempotent algebra with $\idrk(A)=i.$

Here $Z(A)$ is spanned over $k$ by
$\{e_{mn}\mid 2\leq m,n\leq i\},$ hence is nonzero if $i\geq 2.$
\end{lemma}

\begin{proof}
It is easy to verify that for any element $e$ of any associative
algebra $A_0,$ the operation $a*b=a\,e\,b$ is again associative.
(When $e$ is not a right zero divisor, the resulting algebra
$A$ is isomorphic as an algebra to the right ideal
$A_0 e\subseteq A_0 ,$ via the map $a\mapsto ae.)$
So our $A$ is an associative $\!k\!$-algebra.
To verify idempotence, we note that each basis element $e_{mn}$ is
a product, $e_{m1}*e_{1n}.$

Now recall that the rank, as a
matrix, of any product matrix $S\,T$ in $M_i(k)$
is less than or equal to each of $\mathrm{rank}(S)$
and $\mathrm{rank}(T).$
Hence a product $a*b=a\,e_{11}\,b$ has rank
$\leq\mathrm{rank}(e_{11})=1$ as a matrix.
Thus, to get a matrix of rank~$i,$ such as the
identity matrix, we need at least $i$ summands.

To show that $i$ summands always suffice, recall that every
matrix of rank $1$ can be written $uv$ for some column matrix
$u$ and row matrix $v.$
If we embed $u$ as the first column of a matrix $u'\in A,$ and $v$ as
the first row of some $v'\in A,$ we see that $u'*v'=uv.$
Since every $i\times i$ matrix $w$ is a sum of at most $i$ rank-$\!1\!$
matrices (e.g., the $i$ matrices that agree in one column with $w,$
and have zeroes everywhere else),
every matrix $w$ is the sum of $i$ products in $A.$

It follows immediately from the definition of our
multiplication that the elements $e_{mn}$ with $2\leq m,n\leq i$ are
in $Z(A),$ and it is easy to see that any element not in the span
of these elements
is sent to a nonzero value either by left or right multiplication
in $A$ by $e_{11},$ giving the asserted description of $Z(A).$
\end{proof}

Here is the closely related Lie example, though we will not attempt to
determine its idempotence rank and total annihilator ideal as precisely
as in the above case.

\begin{lemma}\label{L.ae1122b}
For any field $k$ of characteristic not $2,$ and any integer $i\geq 2,$
let $A$ be the $\!k\!$-algebra with underlying $\!k\!$-vector-space
the subspace of $M_i(k)$ consisting of all matrices
in which the coefficients of $e_{11}$ and of $e_{22}$ sum to zero,
and with operation given by
\begin{equation}\begin{minipage}[c]{35pc}\label{d.e11e22}
$[a,b]\,=\,a\,(e_{11}+e_{22})\,b-b\,(e_{11}+e_{22})\,a.$
\end{minipage}\end{equation}

Then $A$ is an idempotent Lie algebra with
$i/4\,\leq\,\idrk(A)\,\leq\,i+3.$
\textup{(}Thus, taking $i$ sufficiently large,
we get arbitrarily large idempotence ranks.\textup{)}

$Z(A)$ contains all elements
$\{e_{mn}\mid 3\leq m,n\leq i\},$ hence is nonzero if $i\geq 3.$
\end{lemma}

\begin{proof}
By the general observation with which we began
the proof of the preceding lemma,
$M_i(k)$ is an associative algebra under
the multiplication $a*b=a\,(e_{11}+e_{22})\,b;$ hence it becomes a
Lie algebra under the corresponding commutator
bracket operation~(\ref{d.e11e22}).

It is easy to check that the set of basis elements $e_{mn}$ in which
one or both of $m,\ n$ are $\geq 3$ spans a $\!2\!$-sided ideal
in the above associative algebra structure.
(In a product $a*b,$ the only way such a basis element
occurring in $a$ or $b$ can lead to a
nonzero term of the product is when any index $\geq 3$
is ``facing away from'' the
factor $e_{11}+e_{22}$ in the definition of our multiplication; hence
such an index survives in every term of the product.)
Consequently, in examining the coefficients of $e_{11}$ and $e_{22}$
in a commutator $[a,b]=a*b-b*a,$ we can without loss of generality
assume that all elements $e_{mn}$ occurring in the
expressions for $a$ and $b$ have both subscripts in $\{1,\,2\}.$
Thus, we are reduced to computing in $M_2(k),$ and there our
multiplication is the ordinary multiplication, hence our brackets
are ordinary commutator brackets, and we know that the value of any such
bracket has trace zero.
So the range of our bracket operation on $M_i(k)$ contains only
matrices in which the coefficients of $e_{11}$ and $e_{22}$
sum to zero, hence the set of matrices with that
property indeed forms a Lie algebra $A.$

This Lie algebra contains the simple, hence
idempotent, Lie subalgebra $\mathrm{sl}_2(k),$ so to show
$A$ is idempotent, it will suffice to show that
the range of the bracket also contains all matrix units $e_{mn}$
with at least one of $m,\ n\geq 3.$
For $n\geq 3,$ we have
\begin{equation}\begin{minipage}[c]{35pc}\label{d.e1n}
$e_{1n}\,=\,[e_{11}\,{-}\,e_{22},\ e_{1n}].$
\end{minipage}\end{equation}
Elements $e_{2n},$ $e_{m1}$ and $e_{m2}$ are obtained by obvious
variants of this calculation,
while if both $m$ and $n$ are $\geq 3,$ we have
\begin{equation}\begin{minipage}[c]{35pc}\label{d.emn}
$e_{mn}\,=\,[e_{m1},\,e_{1n}],$
\end{minipage}\end{equation}
completing the proof of idempotence.

Every bracket $[a,b]$ is by
definition a difference of two matrices each of which,
under the ordinary matrix multiplication of $M_i(i),$ has an internal
factor $e_{11}+e_{22},$ hence both of which have rank $\leq 2.$
Thus, $[a,b]$ has rank $\leq 4,$ so at least $i/4$ summands
are needed to get an element of rank~$i.$
This gives our lower bound on $\idrk(A).$

To get the upper bound, note first that if in~(\ref{d.emn}) we hold
$n$ fixed, and taken an arbitrary
$\!k\!$-linear combination of the resulting equations
for all $m\geq 3,$ we get, as a single
commutator $[x,\,e_{1n}],$ an arbitrary
column in position $n\geq 3$ having top two components zero.
Thus, summing $i-2$ such commutators, we can get any matrix living
in the lower right-hand $i-2\times i-2$ block of $M_i(k).$
Linear combinations of the equations~(\ref{d.e1n}), and of the
three variants mentioned following it, show that
with four more commutators, we can fill in everything but the
upper left $2\times 2$ block of a general member of $A.$
Since every element of $\mathrm{sl}_2(k)$ is a commutator,
we can fill in that block in one more step; so every member
of $A$ is a sum of $(i-2)+4+1=i+3$ commutators.

The final sentence of the lemma follows from the observation that
the elements $e_{mn}$ with $3\leq m,n\leq i$
lie in the total annihilator ideal
of the associative multiplication $a(e_{11}+e_{22})b,$
hence a fortiori in the total annihilator ideal of our Lie bracket.
\end{proof}

The final example of this group,
giving nonassociative, non-Lie, finite-dimensional {\em simple}
algebras of unbounded idempotence rank, plays further changes on
the idea of a multiplication whose outputs are rank-$\!1\!$ matrices.

\begin{lemma}\label{L.knxMn}
For any field $k$ and positive integer $i,$ let $A$ be
the $\!k\!$-algebra with underlying $\!k\!$-vector-space
$k^i\times M_i(k),$ and with multiplication defined as follows:\\
-- For $u,\,v\in k^i,$ written as row vectors, $u*v$ is the matrix
$u^\mathrm{T}\,v\in M_i(k),$ where $^\mathrm{T}$ denotes transpose.\\
-- For $S,\,T\in M_i(k),$ $S*T$ is the vector in $k^i$ whose
$\!m\!$-th entry is the $\!(m-1)\!$-st main-diagonal
entry of the ordinary matrix product $S\,T.$
Here we treat subscripts cyclically, so that for $m=1,$
the ``$\!(m-1)\!$-st'' main-diagonal entry means the $\!i\!$-th.\\
-- Products, in either order, of a member of $k^i$ and a member
of $M_i(k)$ are zero.

Then $\idrk(A)=i,$ and $A$ is simple.
\end{lemma}

\begin{proof}
In the product of any two members of $A,$ the $\!M_i(k)\!$-component
will be a matrix $u^\mathrm{T}\,v$ $(u,v\in k^i),$
hence will have rank $\leq 1;$
so at least $i$ such products must be summed to get elements whose
$\!M_i(k)\!$-components have rank $i.$
To get an arbitrary element $(v, S)$ as a sum
$a_1*b_1+\dots+a_i*b_i,$ one first selects, as the $\!k^i\!$-components
of $a_1,b_1,\dots,a_i,b_i,$ pairs of vectors having
products, under our multiplication, summing to $S.$
In all but one of these pairs, one then takes the
$\!M_i(k)\!$-components zero, and in the remaining pair, one takes for
those components matrices $S$ and $T$ such that $S*T$ is
the desired first component $v.$
Thus, $\idrk(A)=i.$

To show that $A$ is simple, let $C$ be a nonzero ideal.
Then $C$ either contains an element with nonzero
$\!k^i\!$-component, or an element with nonzero $\!M_i(k)\!$-component.

In the former case, the square of the element in question will
have nonzero $\!M_i(k)\!$-component, so in either case
$C$ contains an element of the latter sort; say $(v,S).$
Suppose the matrix $S$ has nonzero $(m,m')$
entry, which we may assume without loss of generality is $1.$

Let us form the product $(v,S)*(0,e_{m'm}).$
This will have $\!M_i(k)\!$-component zero;
to determine its $\!k^i\!$-component, note that
the only nonzero main-diagonal component of
$S\,e_{m'm}$ is a $1$ in the $\!m\!$-th position.
So by our description of products $S*T,$ the element
$S*e_{m'm}\in k^i$ will be the vector $f_{m+1}$
with a $1$ in the $\!(m+1)\!$-st position and zeroes elsewhere.
So $C$ contains $(f_{m+1},0).$
Squaring this, we get $(0,e_{m+1,m+1}),$ and squaring
that in turn gives $(f_{m+2},0).$
Repeating this process $i$ times
(and recalling that in this computation, subscripts are
treated cyclically), we get all of $(f_1,0),\dots,(f_i,0).$

Taking linear combinations of these gives all elements $(v,0);$
multiplying pairs of such elements, and adding together
families of $i$ such products, gives all elements $(0,S);$
adding these two sorts of elements
we get all of $A,$ completing the proof of simplicity.
\end{proof}

\subsection{Idempotence rank and base change.}\label{S2.idrk_neq}

We will now give an example showing that the idempotence rank
of a finite-dimensional algebra can
go down (or up) under extension of base field.
However, our example will be non-simple and non-Lie.

We arrived at this example by looking, first, for an
example of this sort for the invariant $\mbox{max-rank}(m)$
of a map $m:A\otimes_k B\to C$ of vector spaces over a field $k,$
defined in~(\ref{d.maxrank}),
a generalization of the idempotence rank of an algebra.
We wondered whether we could find such a map $m$
for a general field $k,$ which, when restricted
to tensors of rank $\leq 1,$ would be surjective if $k$ was
the complex numbers, but such that over the real numbers, the
range would be constrained by inequalities in the coordinates; so that
on passing from the reals to the complexes, the value
of $\mbox{max-rank}(m)$ would drop from a larger value to~$1.$
A little thought shows that for this to happen,
$A$ and $B$ must each be at least $\!2\!$-dimensional,
so we tried $A=B=k^2.$
The tensors of rank $\leq 1$ within $k^2\otimes_k k^2$ can be pictured
as the matrices of rank $\leq 1$ in $M_2(k),$
a set with three degrees of freedom, suggesting that a linear image
of this set in $k^3$ might have the desired properties.

It turns out that if we map a $2\times 2$ matrix $((a_{mn}))$
to the $\!3\!$-tuple consisting of its upper-right and lower-left
entries, and its trace, this has the desired properties.
Indeed, for a $\!3\!$-tuple $(a_{12},a_{21},t)$ of
elements of $k$ to arise in this way from a matrix $((a_{mn}))$ of
rank $\leq 1,$ the entries $a_{11}$ and $a_{22}$ of the latter
matrix must have product $a_{12}a_{21}$ (to make the determinant
$a_{11}a_{22}-a_{12}a_{21}$ zero)
and must sum to $t$ (by definition of the trace).
But two elements of $k$ having sum $t$ and
product $a_{12}a_{21}$ must be the
roots of the quadratic polynomial $x^2-tx+a_{12}a_{21}.$
Over the complexes, such a polynomial always has roots, but it will
have roots over the reals only when the
inequality $t^2-4\,a_{12}a_{21}\geq 0$ holds.

To embody this idea in an idempotent algebra,
let $k$ be any field of characteristic not $2,$
and $A$ the $\!k\!$-algebra with
underlying vector space $k^3,$ and multiplication
\begin{equation}\begin{minipage}[c]{35pc}\label{d.rkchg}
$(a,b,c)*(a',b',c')\,=\,(ab',\ ba',\ aa'+bb').$
\end{minipage}\end{equation}
Note that the components of the product are the
$(1,2)$ entry, the $(2,1)$ entry, and the trace, of
the $2\times 2$ rank-$\!1\!$ matrix $(a,b)^\mathrm{T}\,(a',b').$
It thus follows from the preceding observations
that an element $(r,s,t)\in A$ is a product in $A$
if and only if $t^2-4\,rs$ is a square in $k.$
Hence, if $k$ is algebraically closed (or even if it
is quadratically closed, i.e., if every element of $k$ is a square), we
get $\idrk(A)=1.$

Conversely, we see
that if $k$ is not quadratically closed, $\idrk(A)$ will not be $1.$
Rather, it turns out that it
is $2,$ since any element $(r,s,t)$ can be written
$(r,0,t)+(0,s,0),$ and each of these summands
satisfies our criterion for being a product in $A.$
In particular, if we construct the above algebra over the field
of reals, and then extend scalars to the complexes, the
idempotence rank drops from $2$ to~$1.$

If, inversely, we start with a quadratically closed
field $k,$ and extend scalars to a non-quadratically closed
field $K\supseteq k,$ then $\idrk(A)$ will increase from~$1$ to~$2.$
(If $k$ is algebraically closed, we must, of course, take $K$
transcendental to get a proper extension.
However, a general {\em quadratically}
closed field $k$ can have {\em finite algebraic} extensions $K$ which
are not quadratically closed \cite[Corollary 7.11(1)]{lam}.)

(Incidentally, our use of the term ``quadratically closed'',
defined above, follows \cite{lam}, but is distinct from the usage
in \cite[p.462, Exercises 8-9]{SL.Alg}, where it means
that for each $c\in K,$ one of $c$ or $-c$ is a square.
Evidently, one definition is modeled on the properties of a
subfield of $\C$ closed under taking square roots,
the other on the properties of a square-root-closed
subfield of $\R.)$

\subsection{Lie examples based on inseparability.}\label{S2.non_alg_cl}

To build examples of simple Lie algebras in positive characteristic
which misbehave under change of base field, let us start with
examples that don't misbehave.
Namely,
\begin{equation}\begin{minipage}[c]{35pc}\label{d.properly}
Let $k_0$ be a field of characteristic $p>0,$ and $L_0$ a
finite-dimensional simple Lie algebra over $k_0,$ such that
for every extension field $k$ of $k_0,$ the Lie algebra
$L_0\otimes_{k_0} k$ is again simple.
\end{minipage}\end{equation}

For instance, we can take $L_0=\mathrm{sl}_n(k_0)$ for any $n\geq 2$
relatively prime to $p.$

Recall now that if $K$ is a finite {\em inseparable}
extension of a field
$k,$ then the commutative ring $K\otimes_k K$ has nilpotent elements.
(This is easily seen if $K$ is purely inseparable, and so has an
element $x\notin k$ with $x^p\in k:$
then $(x\otimes 1 - 1\otimes x)^p = x^p-x^p=0.$
To see the general case, recall \cite[Proposition~V.6.6, p.250]{SL.Alg}
that there will exist a separable subextension $F\subseteq K$
over which $K$ is purely inseparable.
By the preceding observation, $K\otimes_F K$ has
nilpotent elements; but that ring is a homomorphic image
of $K\otimes_k K,$ and if a homomorphic image of a finite-dimensional
algebra over a field has nilpotents, the
original algebra must also have them.)
From this we get

\begin{lemma}\label{L.bad_extn}
For $k_0$ and $L_0$ as in~\textup{(\ref{d.properly})}, let
$k\subseteq K$ be
extension fields of $k_0,$ with $K$ a finite inseparable extension
of $k,$ and let $L=L_0\otimes_{k_0} K,$ regarded as an
algebra over $k\subseteq K.$

Then $L$ is a finite-dimensional simple Lie algebra over $k,$ but
$L\otimes_k K$ is not semisimple \textup{(}i.e., it has
a nonzero nilpotent ideal\textup{).}
\end{lemma}

\begin{proof}
$L$ is clearly a finite-dimensional Lie algebra over $k,$
and is simple by~(\ref{d.properly}).

Note that $L\otimes_k K=(L_0\otimes_{k_0} K)\otimes_k K\cong
L_0\otimes_{k_0}(K\otimes_k K).$
Since $K$ is inseparable over $k,$ we can
find a nonzero nilpotent element $\varepsilon\in K\otimes_k K.$
Thus, $L_0\otimes \varepsilon K$ is a nilpotent ideal
of $L_0\otimes_{k_0}(K\otimes_k K).$
\end{proof}

Recall next that in the proof of Theorem~\ref{T.x1x2},
we were able to pull the property $L=[x_1,L]+[x_2,L]$ down from
the case of a field $K$ to that of an infinite subfield $k.$
The next lemma shows that the condition of being generated as an
algebra by two elements, from which we proved that property,
cannot be pulled down in that way.

Though as is well-known, any finite separable extension field $K$
of a field $k$ can be generated over $k$ by a single element
(the Theorem of the Primitive Element
\cite[Theorem~V.6.6, p.243]{SL.Alg}), we shall use in our
construction the fact this fails arbitrarily
badly for inseparable extensions.
For instance, if we take for $k$ a pure transcendental
extension $k_0(t_1,\dots,t_N)$ of $k_0,$
then the degree-$\!p^N\!$ extension
$K=k(t_1^{1/p},\dots,t_N^{1/p})$ cannot be generated
over $k$ by fewer than $N$ elements
\cite[Theorem~8.6.4]{GDVS}.

\begin{lemma}\label{L.not-2-gen}
Given $k_0$ and $L_0$ as in~\textup{(\ref{d.properly})}, let
$d=\dim_{k_0}(L_0),$ let $n$ be any positive integer, and
let $k\subseteq K$ be extension fields of $k_0,$
such that $K$ is finite over $k,$ but cannot
be generated over $k$ by fewer than $nd+1$ elements.
Again, let $L=L_0\otimes_{k_0} K,$ regarded
as a finite-dimensional simple Lie algebra over $k\subseteq K.$

Then $L$ cannot be generated as a Lie algebra over~$k$ by
fewer than $n+1$ elements.
\end{lemma}

\begin{proof}
Let $B=\{b_1,\dots,b_d\}$ be a $\!k_0\!$-basis for $L_0.$
Then $B$ will likewise be a $\!K\!$-basis for $L.$
Given $n$ elements $x_1,\dots,x_n\in L,$ their expressions
in terms of that basis will involve $d\,n$ coefficients in $K.$
By assumption, $d\,n$ elements cannot generate $K$ over $k,$
so those coefficients lie in a proper subextension $F\subseteq K.$
Thus $x_1,\dots,x_n$ lie in $L_0\otimes_{k_0} F,$ a
proper $\!k\!$-subalgebra of $L_0\otimes_{k_0} K=L.$
\end{proof}

\section{Some questions, and some directions for further study.}\label{S.questions}
\subsection{Can our cardinality restrictions be weakened?}\label{S2.card}

In the main results of this paper, we have assumed
the field $k$ infinite.
Some of those results remain formally true -- but become trivial --
for finite $k:$ the hypothesis that $\card(I)$
be less than any measurable cardinal $>\card(k)$ then says that
$I$ is finite.
(Recall that under the definition we are
following, $\aleph_0$ is a measurable cardinal.)

In Theorem~\ref{T.solv} and Corollary~\ref{C.solv},
we assumed, slightly more generally, that {\em either}
$\card(I)$ {\em or} the supremum of the
dimensions of all the $A_i$ was less than all measurable
cardinals greater than $\card(k).$
What this would say for finite $k$
is that either $\card(I)$ is finite, or the
dimensions of the $A_i$ have a common finite bound.
Under this assumption, the conclusions of those two results
follow easily from Proposition~\ref{P.via_ultra} and
the well-known fact that any ultraproduct of
{\em finite} algebraic structures (with only finitely many operations),
of bounded cardinalities, is isomorphic to one of those structures.

What we would like to know, of course, is

\begin{question}\label{Q.finite}
Suppose $k$ is a finite field, and $f:\prod_I A_i\to B$ a homomorphism
of $\!k\!$-algebras, with $I$ infinite, and no common finite
bound assumed on the $\!k\!$-dimensions of the $A_i.$
Do some or all of the main results of this paper
\textup{(}other than Theorem~\ref{T.nilp}\textup{)}
have versions valid for this case?
\textup{(}Or can some other results in the same spirit be
established?\textup{)}
\end{question}

We have excluded Theorem~\ref{T.nilp} because, as mentioned,
a result on nilpotent algebras with no condition that $k$ be infinite
is indeed proved in~\cite{pro-np}.
(Using that result, and the close relationship between
nilpotence and solvability for finite-dimensional Lie algebras
in characteristic~$0,$ an analog of Theorem~\ref{T.solv2} for that
particular case is also obtained in~\cite{pro-np}.)

Another sort of size restriction in our results on
homomorphisms $\prod_I A_i\to B$ concerned the object $B.$
Here {\em some} restriction is needed, since if we allowed $B$ to
be $\prod_I A_i,$ the identity map of that algebra would be a
counterexample to most of our results.
But it is not clear that the conditions need to
be as strong as those we have used.
For instance

\begin{question}\label{Q.no_CC}
Does Theorem~\ref{T.via_A_i} remain true if we delete the
hypothesis that $B$ satisfy chain condition on almost direct factors?
\end{question}

In \cite[Theorem~9(i-ii)]{prod_Lie2} we indeed prove results like
Theorem~\ref{T.via_A_i} without the chain
condition hypothesis -- but having,
instead, restrictions such as $\card(I)\leq\card(k).$
So we want to know whether we can do without either sort of condition.

A result in which the condition on the codomain algebra might
be weakened in a different way
is Theorem~\ref{T.solv} above, where the codomain
is assumed both simple (much stronger than having
chain condition on almost direct factors) and countable-dimensional,
and we do not know whether the latter condition can be dropped.
We also don't know, in that case, whether we need the restriction on
the size of the algebras $A_i$ or the
index-set $I$ relative to measurable cardinals.
So we ask

\begin{question}\label{Q.solv}
If $(A_i)_{i\in I}$ is a family of algebras over an infinite field
$k,$ such that no $A_i$ admits a homomorphism onto a simple algebra,
can $\prod_I A_i$ admit a homomorphism onto a simple algebra?
\end{question}

Let us note that in the existing proof of Theorem~\ref{T.solv}, the
dimension-restriction on the codomain {\em can} be slightly weakened:
Using the full strength of Theorem~\ref{T.bigger}, we see that if
$k$ has cardinality $>\aleph_1,$ we get the indicated nonexistence
result (with the parenthetical generalization in the first
sentence appropriately adjusted)
not just for homomorphisms onto simple algebras that are
countable-dimensional,
but onto the larger class of simple
algebras of $\!k\!$-dimension $<\card(k).$

Concerning the hypothesis in that result that the dimensions
of the $A_i$ be less than any measurable cardinal $>\card(k),$
we wonder whether one might be able
to remove this by showing that simple algebras
$A_i$ of such large dimensions can be replaced by
simple subalgebras of smaller dimensions, without affecting the
desired properties.

In \cite[\S8]{pro-np}, examples are given of homomorphic images $B$ of
inverse limits $A$ of nilpotent algebras in which
$B$ has various properties that inverse limits of
nilpotent algebras cannot themselves have;
e.g., an associative example where $B$ contains
a nonzero element $y$ such that $y\in ByB,$
and a nonassociative example where $B$ contains
an element such that $y^2=y.$
However, the analog of Question~\ref{Q.solv}
for inverse limits of nilpotent is open;
it is part of \cite[Question~23]{pro-np}.

\subsection{On idempotence rank.}\label{S2.1bracket}
The next question poses the problem that we
skirted by obtaining our bound on idempotence ranks
of simple Lie algebras using~\cite{Bois} instead of~\cite{Brown}.

\begin{question}[{also asked in~\cite[pp.652-653]{KHH+SAM} for $k=\R$}]\label{Q.1bracket}
Does every finite-dimensional simple Lie
algebra $L$ over an infinite field $k$ have idempotence rank~$1$?
\end{question}

\subsection{On the chain condition on almost direct factors.}\label{S2.CCq}
In~\S\ref{S2.CC} we saw that a finitely generated associative
algebra over a
field need not have chain condition on almost direct factors.
On the other hand, a finitely generated {\em commutative} associative
algebra over a field is Noetherian, and so does
have that chain condition.

\begin{question}\label{Q.CC}
Does every finitely generated {\em Lie}
algebra over a field have chain condition on almost direct factors?
\end{question}

\subsection{Idempotence rank, number of generators, and base change.}\label{S2.idrkq}
The algebras of~\S\ref{S2.idrk_neq},
whose idempotence ranks could increase and decrease under base change,
were neither associative nor Lie.
Also, the motivating idea of that example -- an image-set which,
when the base field is $\R,$ is constrained by inequalities,
but which is not so constrained when the base field is $\C$
-- only seems to lead to examples where
the idempotence rank changes by~$1.$
So we ask

\begin{question}\label{Q.idp_rk_neq}
\textup{(a)}\ \ Do there exist finite-dimensional {\em associative or
Lie} algebras over an infinite field
whose idempotence ranks change under base extension?

\textup{(b)}\ \ Do there exist finite-dimensional algebras of any sort
over an infinite field
whose idempotence ranks change by more than~$1$ under base extension?
\end{question}

In \cite[\S2]{x1Adots}, examples are given of finite-dimensional
(nonassociative, non-Lie)
algebras over {\em finite} fields whose idempotence ranks change
by arbitrarily large amounts under base
extension, and of an infinite-dimensional
{\em commutative associative} algebra over $\R$
whose idempotence rank goes down (though only by~$1)$
on extension of scalars to $\C.$

In the same vein is

\begin{question}[J.-M.\,Bois, personal communication]\label{Q.2-gen:finite}
Let $k$ be a finite field, $K$ its algebraic closure, and
$L$ a finite-dimensional Lie algebra over $k.$
If $L\otimes_k K$ can be generated over $K$ by two elements,
can $L$ be generated over $k$ by two elements?

For instance, is $\mathrm{sl}_n(k)$ generated over $k$ by
two elements for all $n\geq 2$ relatively prime to $\mathrm{char}(k)$?
\end{question}

\subsection{Centroids to the rescue for our inseparable-extension examples?}\label{S2.centroid}

In \S\ref{S2.non_alg_cl}, where we constructed simple
Lie algebras in positive characteristic that ``misbehaved''
under base change, the trick was to
treat them as having base field $k,$ though they were Lie
algebras over a larger field $K,$ which was inseparable over $k.$
As $\!K\!$-algebras, they are well-behaved.

If $L$ is a finite-dimensional simple Lie algebra over a field $k,$
the largest field $K$ to which the Lie structure extends,
called the {\em centroid} of $L,$ consists of
the $\!k\!$-linear endomorphisms $\varphi$
of $L$ which respect all the adjoint maps $[x,-\,]$ $(x\in L),$
i.e., which satisfy $\varphi([x,y])=[x,\varphi(y)]$ for
all $x,y\in L$ \cite[p.290]{NJ}.
If $K=k,$ $L$ is called a {\em central} simple Lie algebra.
Thus, every finite-dimensional simple Lie algebra over a field is a
central simple Lie algebra over its centroid.
It is plausible that if we of look at our Lie
algebras as algebras over their centroids, the kind of misbehavior
obtained in~\S\ref{S2.non_alg_cl} will not occur:

\begin{question}\label{Q.centroid}
Let $L$ be a finite-dimensional {\em central} simple Lie algebra over
a field $k$ of characteristic not $2$ or~$3.$

\textup{(a)}\ \ Can $L$ be generated over $k$ by two elements?

\textup{(b)}\ \ Will $L\otimes_k K$ be a direct product of
simple Lie algebras over $K$ for all extension fields $K$ of $k$?

\textup{(c)}\ \ If either of the above questions has a positive
answer, does it remain so if rather than assuming $L$ central,
we merely assume the centroid of $L$ to be separable over $k$?
\end{question}

\subsection{Characteristics $2$ and $3$}\label{S2.char=2,3}
Of course, we would like to know

\begin{question}\label{Q.char2,3}
Does the conclusion of Theorem~\ref{T.simple_Lie}\textup{(i)}
hold in the excluded characteristics,~$2$ and~$3$?
\end{question}

Perhaps, when the structure theory is
extended to those last two characteristics,
the result of \cite{Bois} that we used in the proof
will also go over, yielding an affirmative answer.
On the other hand, a weaker result than that of \cite{Bois},
perhaps asserting generation by $3$ or $4$ elements rather
than~$2,$ might be easier to prove than the optimal result,
and might not require a full structure theory.
\vspace{.5em}

We give the last three points of this section as topics
to be investigated, rather than formal questions.

\subsection{Variant formulations of solvability and nilpotence.}\label{S2.solv_&_nil}
Of the two versions of our result on homomorphic images of direct
products of solvable Lie algebras, we got the first,
Corollary~\ref{C.solv}, using the criterion that a finite-dimensional
Lie algebra over a field of characteristic~$0$
is solvable if and only if it admits no homomorphism
onto a simple Lie algebra, while for the second,
Theorem~\ref{T.solv2}, we used the characterization of solvability
(in any characteristic) by a disjunction of identities.
Thus, our results on Lie algebras were obtained as cases
of two different results on general algebras.

There are other elegant characterizations of solvability of a
finite-dimensional Lie algebra $L:$ in arbitrary characteristic,
the condition that $L$ have no nontrivial idempotent subalgebra;
in characteristic~$0,$ either the
condition that the ideal $[L,L]$ be nilpotent
(which is sufficient but not necessary in general characteristic),
or that $L$ contain no simple subalgebra (necessary, but not
sufficient in general characteristic).

Likewise, finite-dimensional {\em nilpotent} Lie algebras $L$
can be characterized among
finite-dimensional Lie algebras in other ways than
the one used in Theorem~\ref{T.nilp}: as those
with no nonzero ideals $C$ such that $[L,C]=C,$ as
those with no nonzero elements $x$ such that $x$ belongs to the
ideal generated by $[L,\,x],$ and as those whose nonzero
homomorphic images $M$ all have $Z(M)\neq\{0\}.$

It might be of interest to examine how some of these conditions
on general algebras behave under homomorphic images of direct products.

\subsection{Semisimple Lie algebras in positive characteristic.}\label{S2.semisimple}
A Lie algebra $L$ is called {\em semisimple} if it has no
nonzero abelian ideal.
If the base field $k$ has characteristic zero, the
finite-dimensional semisimple Lie algebras are just the
finite direct products of simple Lie algebras,
so our results on homomorphic images of direct products
of simple Lie algebras imply the corresponding statements
for products of semisimple Lie algebras.

When the base field has positive characteristic, a
finite-dimensional semisimple Lie algebra need not be a direct
product of simple Lie algebras \cite[p.133, top paragraph]{Strade}.
The present authors know nothing about their structure.

In particular, what Lie algebras are homomorphic images of
finite-dimensional semisimple Lie algebras?
It is conceivable that all are.
(By analogy, every finite group $G$
is indeed a homomorphic image of a finite group having no abelian
normal subgroup, namely, a wreath product of $G$ with a
finite simple group.)
If so, then little can be said about homomorphic
images of infinite products of such algebras
-- though something might be said about homomorphisms from
infinite products {\em onto} semisimple Lie algebras.

These seem to be questions for the expert in
Lie algebras over fields of positive characteristic.
An introduction to the subject is~\cite{Str+Fa}.

\subsection{Restricted Lie algebras, and other algebras with additional structure.}\label{S2.p-Lie}

For $k$ a field of positive characteristic $p,$
a {\em restricted} Lie algebra or {\em $\!p$-Lie algebra} over $k$
is a Lie algebra given with
an additional operation, $x\mapsto x^{(p)},$ satisfying certain
identities which, in associative $\!k\!$-algebras,
relate the $\!p$-th power map with the
$\!k\!$-module structure and commutator brackets \cite[\S5.7]{NJ}.
(The concept can be motivated by the
observation that in characteristic $p,$
the set of derivations of an algebra $A$ is closed, in the associative
algebra of $\!k\!$-vector-space endomorphisms of $A,$ not only
under the vector space operations and
commutator brackets, but also under taking $\!p$-th powers.)

Thus, $\!p$-Lie algebras are
not algebras as we define the term in~\S\ref{S.intro}.
Of course, they have Lie algebra structures,
and one can apply our results to those structures.
But a $\!p\!$-Lie algebra may, for instance, be simple under
its $\!p$-Lie algebra structure without being simple under
its ordinary Lie algebra structure.
We leave to others the investigation of homomorphic images of
infinite products of $\!p$-Lie algebras, and, generally, of algebras
with additional operations.

\section{Appendix: Review of ultrafilters and ultraproducts.}\label{S.ultra_defs}
We recall here some standard definitions and notation
(e.g., cf.~\cite[p.211{\em ff}\/]{Ch+Keis}).

A {\em filter} on a nonempty set $I$ means
a family $\F$ of subsets of $I$ such that
\begin{equation}\begin{minipage}[c]{35pc}\label{d.ultra}
$J_1\supseteq J_2\in \F\implies J_1\in \F,$\\[.2em]
\strt$J_1,\ J_2\in \F\implies J_1\cap J_2\in \F.$
\end{minipage}\end{equation}

In view of the first condition,
a filter $\F$ is {\em proper} (not the set of all subsets
of $I)$ if and only if $\emptyset\notin\F.$

A maximal proper filter is called an {\em ultrafilter}\/;
by Zorn's Lemma, every proper filter is contained in an ultrafilter.
It is easy to verify that a proper filter $\U$ is an ultrafilter
if and only if for every $J\subseteq I,$ either $J\in\U$ or $I-J\in\U.$
If $\F$ is a filter (in particular, if it is an ultrafilter) on $I,$ one
says that a subset $J\subseteq I$ is {\em $\!\F\!$-large} if $J\in\F.$
This does not save much ink, but
does help with the intuition of the subject.

If $(A_i)_{i\in I}$ is a family of nonempty sets,
and $\F$ is a filter on the index set $I,$
then the {\em reduced product} $\prod_I A_i\,/\,\F$ is the factor-set
of $\prod_I A_i$ by the equivalence relation that identifies
elements $(a_i)$ and $(a'_i)$ if $\{i\mid a_i=a'_i\}$ is $\!\F\!$-large.
If all the $A_i$ are furnished with operations making them
groups, $\!k\!$-algebras, etc., then this
structure can be seen to carry over to their reduced product,
making the natural map $\prod_I A_i\to\prod_I A_i\,/\F$ a homomorphism.
For any $J\in\F,$ it is not hard to see that
$\prod_I A_i\,/\F\cong\prod_J A_i\,/\F_J,$ where
$\F_J=\{J'\in\F\mid J'\subseteq J\},$ a filter on $J.$
(This observation depends on
our assumption that all $A_i$ are nonempty.)

A reduced product of objects $A_i$ with respect to an
{\em ultrafilter} is called an {\em ultraproduct} of the $A_i.$
An ultraproduct $A^I\,/\,\U$ of copies of a single object $A$
is called an {\em ultrapower} of~$A.$
Note that in this situation, the diagonal image of
$A$ in $A^I$ maps to an isomorphic copy of $A$ within $A^I/\,\U.$

It is known that an ultraproduct $\prod_I A_i\,/\,\U$ satisfies every
first order sentence $s$ which holds on a $\!\U\!$-large subfamily of
the $A_i;$ i.e., for which $\{i\in I\mid A_i$~satisfies $s\}\in\U$
\cite[Theorem~4.1.9(iii)]{Ch+Keis}.

For every $i_0\in I,$ the filter of all subsets of $I$ containing
$i_0$ is called the {\em principal} ultrafilter determined by $i_0.$
The ultraproduct of the $A_i$ with respect to that
ultrafilter is, up to isomorphism, just $A_{i_0}.$
If $I$ is finite, these are the only ultrafilters on $I;$ if $I$ is
infinite, on the other hand, then the
cofinite subsets of $I$ form a proper filter,
so there are ultrafilters containing this filter,
the {\em nonprincipal} ultrafilters.

We note, for perspective, that filters on $I$ correspond to
ideals in the Boolean algebra of subsets of $I,$
by mapping each filter to the ideal of complements of its members.
The ultrafilters correspond to the maximal ideals, which in this
case are the same as the prime ideals.
More generally, for any family of fields $(k_i)_{i\in I},$ the ideals
of $\prod_I k_i$ correspond to the filters on $I,$
each filter $\F$ yielding the ideal of all elements
having $\!\F\!$-large zero-set; in other words, the kernel
of the map from $\prod_I k_i$ to the reduced
product $\prod_I k_i\,/\,\F.$
Again the ultrafilters correspond to the maximal ideals,
and these coincide with the prime ideals.
(The statements about Boolean algebras are essentially the
cases of the statements about products of fields
in which all $k_i$ are the two-element field.)

\section{Appendix: $\!\kappa\!$-complete and non-$\!\kappa\!$-complete ultrafilters.}\label{S.size_ultra}

\begin{definition}[{\cite[p.227]{Ch+Keis}}]\label{D.*k-cplt}
Let $\kappa$ be an infinite cardinal.
Then an ultrafilter $\U$ on a set $I$
is said to be {\em $\!\kappa\!$-complete}
if it is closed under intersections of families of {\em fewer}
than $\kappa$ members.
An $\!\aleph_1\!$-complete ultrafilter \textup{(}i.e., one
closed under countable intersections\textup{)} is called
{\em countably complete}.

An infinite cardinal $\kappa$ is called {\em measurable} if there
exists a nonprincipal $\!\kappa\!$-complete ultrafilter on $\kappa.$
\end{definition}

(We follow~\cite{Ch+Keis} in this definition.
Many authors, restrict the term
``measurable cardinal'' to the case where $\kappa$ is uncountable,
e.g., \cite[p.177]{Drake}.
We shall, rather, explicitly write ``uncountable measurable cardinal''
when that is intended.)

Note that the definition of an ultrafilter makes it
$\!\aleph_0\!$-complete (closed under finite intersections), so
the weakest completeness condition not automatically satisfied
is countable completeness.
Note also that $\aleph_0$ is, under the above definition, a
measurable cardinal, since there exist nonprincipal ultrafilters on it.

It is known \cite[Proposition~4.2.7]{Ch+Keis} that for
any nonprincipal ultrafilter $\U$ on any index set $I,$ there is a
largest cardinal $\kappa$ such that $\U$ is $\!\kappa\!$-complete,
and that this will be a measurable cardinal.
Moreover, if an uncountable measurable cardinal exists, it must
be ``enormous'' in many respects.
In particular, truncating
set theory to exclude it and all larger
cardinals will leave a smaller set theory
that still satisfies the standard axiom system ZFC; hence,
if ZFC is consistent, so is ZFC together with the statement
that there are no uncountable measurable cardinals, and therefore no
nonprincipal
countably complete ultrafilters \cite[Chapter~6, Corollary~1.8]{Drake}.

If uncountable measurable cardinals $\mu$ do exist, then any set $I$
admitting a nonprincipal ultrafilter $\U$ that is $\!\mu\!$-complete for
such a $\mu$ must itself have cardinality at least $\mu$
\cite[Proposition~4.2.2]{Ch+Keis}.
It follows that every element of $\U$ must likewise have
cardinality at least $\mu.$

Thus, the reader may prefer to assume that there are no uncountable
measurable cardinals, or at least that the products of algebras
he or she is interested in will always
be indexed by sets $I$ of less than any uncountable measurable
cardinality, and so read only the first result below, which
concerns {\em non}-$\!\kappa\!$-{\nolinebreak}complete ultrafilters.
But the subsequent results, about $\!\kappa\!$-complete ultrafilters,
show that even if these occur, things still work out fairly
nicely for our purposes!

Note that since the sets not in an ultrafilter $\U$ are the
complements of the sets in $\U,$ the condition that $\U$
be $\!\kappa\!$-complete is equivalent to saying that if a
family of $<\kappa$ sets $I_\alpha\subseteq I$ has union in $\U,$
then at least one of the $I_\alpha$ lies in $\U.$

For $\kappa$ a cardinal, $\kappa^+$ denotes the
successor of $\kappa,$ so that a $\!\kappa^+\!$-complete
ultrafilter is one closed under $\!\kappa\!$-fold intersections;
equivalently, one which always contains some member of a
$\!\kappa\!$-tuple of sets if it contains their union.
We follow the standard convention that every cardinal $\kappa$ is
the set of all ordinals of cardinality $<\kappa,$ hence is itself a set
of cardinality $\kappa.$
The least infinite cardinal, $\aleph_0,$
looked at as the set of natural numbers, is denoted $\omega.$

As mentioned in the preceding section, ultraproducts preserve
first-order sentences; hence an ultraproduct of fields is
not merely a ring, but a field.
Our first result below says that {\em except} in the case involving
large measurable cardinals, a nonprincipal
ultrapower of an infinite field
$k$ is significantly larger than that field.

\begin{theorem}\label{T.bigger}
Let $k$ be an infinite field, let $\kappa=\card(k),$ let $\U$
be a {\em non-$\!\kappa^+\!$-complete} ultrafilter on a set $I,$
and let $K=k^I/\,\U.$
Then the dimension $[K:k]$ is uncountable, and is at least $\card(k).$
\end{theorem}

\begin{proof}
Since $\U$ is not $\!\kappa^+\!$-complete, we can take
a family of $\kappa$ sets $I_\alpha$ $(\alpha\in\kappa)$ which
are not in $\U,$ but whose union is in $\U.$
Deleting from each the union of those that precede it in our indexing,
we may assume that they are disjoint; and throwing in, as
one more set, the complement of their union
(which is not in $\U$ because their union {\em is} in $\U),$
we may assume the $I_\alpha$ have union $I.$
(Some $I_\alpha$ may be empty.)

We shall prove first that $k^I/\,\U$ is transcendental over $k.$
Thus, letting $t$ be a transcendental element, the
elements $(t-c)^{-1}$ $(c\in k)$ will be $\!k\!$-linearly independent,
proving $[K:k]\geq\card(k).$
If $k$ is uncountable, this makes $[K:k]$ uncountable.
On the other hand, for countable $k,$ we shall show
that $k^I/\,\U$ is uncountable, again implying that
$[K:k]$ is uncountable.

To show $k^I/\,\U$ transcendental over $k,$
write $k$ as $\{c_\alpha\mid \alpha\in\kappa\},$ with the
$c_\alpha$ distinct, and let $t\in k^I/\,\U$ be the image of the
element $c\in k^I$ which has value $c_\alpha$ everywhere on $I_\alpha$
for each $\alpha\in\kappa.$
Any nonzero polynomial $p(x)$ over $k$ has only finitely
many roots in $k,$ hence its value at $c$ is zero only on a finite
union of the $I_\alpha,$ hence not on a member of $\U;$ so
$p(t)\neq 0$ in $k^I/\,\U,$ showing that $t$ is transcendental.

On the other hand, suppose $k$ is countable, so that our hypothesis
on $\U$ is that it is not countably complete.
As above, take a decomposition of $I$ into disjoint
sets $I_n\notin\U$ $(n\in\omega).$
We shall show that given any countable list of elements
of $k^I/\,\U,$ we can find an element not in that list,
proving $k^I/\,\U$ uncountable.

Let the members of our list be the images in $k^I/\,\U$ of elements
$a^{(0)},\ a^{(1)},\ \dots\in k^I.$
Then we can choose an element $c\in k^I$ which disagrees with $a^{(0)}$
at each point of $I_0$ (because $k$ has more than one element),
with both $a^{(0)}$ and $a^{(1)}$ at each point of $I_1$
(because $k$ has more than two elements), and so forth.
We then see that for each $n,$ the set at which $c$ agrees with
$a^{(n)}$ is a subset of $I_0\cup\dots\cup I_{n-1}\notin\U;$ hence the
element of $k^I/\,\U$ that $c$ defines is distinct from each member
of the given countable list, as claimed.
\end{proof}

When $\U$ {\em is} $\!\kappa^+\!$-complete, we have a result
of the opposite sort.

\begin{theorem}\label{T.nobigger}
Let $k$ be an infinite field, let $\kappa=\card(k),$ and let
$\U$ be a $\!\kappa^+\!$-complete ultrafilter on a set~$I.$
Then the ultrapower $k^I/\,\U$ coincides with the natural
isomorphic copy of $k$ therein.

In this situation, if $\mu$ is any cardinal $>\kappa$ such that
$\U$ is $\!\mu\!$-complete, and $(A_i)_{i\in I}$ is a family of
$\!k\!$-algebras whose dimensions have supremum $<\mu,$
then the ultraproduct $\prod_I A_i\,/\,\U$
is isomorphic as a $\!k\!$-algebra to~$A_{i_1}$ for some $i_1\in I,$
and the canonical map $\prod_I A_i\to \prod_I A_i\,/\,\U\cong A_{i_1}$
splits \textup{(}is right invertible\textup{)}.
\end{theorem}

\begin{proof}
The first paragraph follows from the case of the second where all the
$A_i$ are one-dimensional, so it suffices to prove the assertions of
the latter paragraph.
We note that these are immediate if $\U$ is principal, so let
us assume the contrary.
In that case, we may take $\mu$ to be the {\em greatest} cardinal
such that $\U$ is $\!\mu\!$-complete.
Thus, $\mu$ is a measurable cardinal.

We note first that for $\lambda$ any cardinal, a
$\!\lambda\!$-dimensional $\!k\!$-algebra can, up to isomorphism,
be taken to have underlying vector space $\bigoplus_\lambda k,$ and its
algebra structure will then be determined by $\lambda^3$ structure
constants $c_{\alpha\beta\gamma}$ $(\alpha,\beta,\gamma\in\lambda),$
where $c_{\alpha\beta\gamma}\in k$ is the coefficient,
in the product of the $\!\alpha\!$-th and $\!\beta\!$-th basis
elements, of the $\!\gamma\!$-th basis element.
(These are subject to the constraint that for each $\alpha$ and
$\beta,$ there are only finitely many $\gamma$ with
$c_{\alpha\beta\gamma}\neq 0;$
but for counting purposes, this will not matter to us.)
Thus, the number of isomorphism classes of $\!\lambda\!$-dimensional
algebras is $\leq\kappa^{\lambda^3}.$
From the fact that a measurable cardinal $\mu$ is inaccessible
\cite[Theorem~4.2.14(i)]{Ch+Keis}, it follows that for any
$\lambda<\mu,$ the cardinality of the set of isomorphism classes
of $\!k\!$-algebras of dimension $\leq\lambda$ is also $<\mu.$
Thus, under the hypotheses of the second paragraph of
our theorem, the $A_i$ fall into $<\mu$ isomorphism classes.

Hence if we partition $I$ according to the isomorphism class of
$A_i,$ the $\!\mu\!$-completeness of $\U$ implies that the subset
$I_0$ corresponding to some one of these classes belongs to $\U.$
Let us assume for notational convenience that all the $A_i$ with
$i\in I_0$ are equal, and call their common value $A_{(0)}.$
We now define the homomorphism
\begin{equation}\begin{minipage}[c]{35pc}\label{d.*q}
$\psi:\ A_{(0)}\ \to\ \prod_I A_i,$\quad where
$\psi(a)_i=a$ if $i\in I_0,$\quad
$\psi(a)_i=0$ otherwise.
\end{minipage}\end{equation}
Composing this with the canonical homomorphism
\begin{equation}\begin{minipage}[c]{35pc}\label{d.*f}
$\varphi:\ \prod_I A_i\ \to\ \prod_I A_i\,/\,\U,$
\end{minipage}\end{equation}
we clearly get an embedding $\varphi\psi:A_{(0)}\to \prod_I A_i\,/\,\U.$

We claim that this is surjective, and hence an isomorphism.
(This is one direction of~\cite[Proposition~4.2.4]{Ch+Keis},
but quick enough to prove directly.)
For $A_{(0)},$ being $\!<\mu\!$-dimensional,
has cardinality $<\mu,$ hence
given any $a=(a_i)_{i\in I}\in A_{(0)}^{I_0},$ we may partition $I_0$
into $<\mu$ subsets $I_{0,b}=\{i\in I_0\mid a_i=b\}$ $(b\in A_{(0)}).$
Again, by $\!\mu\!$-completeness one of these lies in $\U;$
say $I_{0,b_0}.$
It follows that $a$ falls together with $\psi(b_0)$ under $\varphi,$
proving surjectivity of $\varphi\psi:A_{(0)}\to\prod_I A_i/\U.$
Thus, $\varphi\psi$ is an isomorphism,
so in particular, $\varphi$ splits.
\end{proof}

If one deletes the bound assumed above on the dimensions of the $A_i,$
one can still show that $\prod_I A_i\,/\,\U$ has properties very
close to those algebras, as illustrated by the next two results.
These are instances of \cite[Theorem~4.2.11]{Ch+Keis}, which says
that the fact quoted earlier,
that ultraproducts preserve first-order sentences
satisfied on $\!\U\!$-large sets, can be strengthened, for
$\!\kappa\!$-complete ultrafilters, to refer to an extended first-order
language allowing conjunctions and disjunctions of all families
of fewer than $\kappa$ sentences.
(Incidentally, in the first result below, we could replace
each identity $W_m=0$ with an arbitrary set of identities;
but for simplicity of statement we refer to single identities, the
case needed in in \S\ref{S.nilp+}.)

In the preceding theorem, the assumption that $k$ be a field was
not essential, but a wordier statement and proof would have
been needed without it.
The next two results are
easy to state and prove for general $k,$
so we revert to our default assumption that $k$ is any commutative ring.

\begin{proposition}[{cf.\ \cite[Theorem~4.2.11]{Ch+Keis}}]\label{P.ids}
Let $(A_i)_{i\in I}$ be a family of $\!k\!$-algebras,
and $\U$ a countably complete ultrafilter on the index set $I.$

Suppose $W_1=0,$ $W_2=0,$ $W_3=0,\,\dots$ is a countable list of
identities for $\!k\!$-algebras, such that each $A_i$
satisfies at least one of these identities.
Then the ultraproduct $\prod_I A_i\,/\,\U$
satisfies at least one of those identities.
\end{proposition}

\begin{proof}
For $m=1,\,2,\,\dots\,,$ let $J_m$ be the set of $i\in I$ such
that $A_i$ satisfies the identity $W_m=0.$
By assumption, $I=\bigcup_m J_m,$ hence since $\U$ is
countably complete, there is an $m$ such that $J_m\in\U.$
Hence $\prod_{J_m} A_i$ satisfies $W_m=0,$ hence so does its image,
$\prod_I A_i\,/\,\U.$
\end{proof}

\begin{proposition}\label{P.simple}
Let $(A_i)_{i\in I}$ be a nonempty family of $\!k\!$-algebras,
and $\U$ a countably complete ultrafilter on the index set $I.$
Then if all $A_i$ are simple $\!k\!$-algebras, so
is the ultraproduct $\prod_I A_i\,/\,\U.$
\end{proposition}

\begin{proof}
An algebra $A$ is simple
if and only if~(i)~$A\neq\{0\},$ and~(ii)~for every nonzero
$a\in A,$ the ideal of $A$ generated by the set $aA+Aa$ is all of $A.$
Clearly,~(i) carries over from the $A_i$ to their ultraproduct.
Note that~(ii) means that for every nonzero $a\in A,$ every $b\in A$
can be represented as a sum of products of elements of $A,$
all of length $\geq 2,$ in which each product includes a factor $a.$
(We have made no reference in this statement to coefficients in $k.$
This was our point in using products of length $\geq 2:$
After selecting an instance of $a$ in each product, we can absorb
a coefficient from $k,$ if any, into one of the {\em other} factors.)

Now there are only countably many forms such an expression as a
sum of products can take.
Indeed, to generate a form for such an
expression, one chooses the natural number
that is to be the number of summands, for each summand one
chooses the natural number $\geq 2$ that is to be its length, and
given this length, one chooses one of the finitely many
positions for the factor $a$ to appear in,
and, finally, one of the finitely many ways
for the (nonassociative) product to be bracketed.

Now let us be given $(a_i)\in\prod_I A_i$ with nonzero image in
$\prod_I A_i\,/\,\U,$ and any $(b_i)\in\prod_I A_i.$
Let $J=\{i\in I\mid a_i\neq 0\}.$
Since the image of
$(a_i)$ in our ultraproduct is assumed nonzero, we have $J\in\U.$
For each $i\in J,$ the simplicity of $A_i$ says
that we can write $b_i$ as a sum of products of lengths $\geq 2$ in
elements of $A_i,$ such that each of these products has a factor $a_i.$
Choosing for each $i$ such an expression for $b_i,$ we can now partition
$J$ as $\bigcup_{m\in\omega} J_m$ according to the form of this
expression, since we have seen that there
are only countably many such forms.
By countable completeness, for some $m$ we have $J_m\in\U.$
Suppose our $\!m\!$-th expression involves $n$ variables other than $a.$
Then we can choose $n$ elements of
$\prod_I A_i$ which, at every $i\in J_m,$ represent the values
used in our expression for $b_i.$
The images of these elements in $\prod_I A_i\,/\,\U$ will
therefore witness the condition that the image of $(b_i)$
lies in the ideal generated by the image of $(a_i).$
This proves the simplicity of $\prod_I A_i\,/\,\U.$
\end{proof}

The above result is not true for non-countably-complete ultrafilters.
For instance, if $k$ is any field, then within the product
of matrix rings $\prod_{i\in\omega} M_i(k),$ the set of elements
$(a_i)_{i\in\omega}$ $(a_i\in M_i(k))$ such that the set of
integers $\{\mathrm{rank}(a_i)\mid i\in\omega\}$
is {\em bounded} forms an ideal.
(Closure under multiplication by arbitrary elements of $\prod M_i(k)$
is immediate; closure under addition
follows from the observation that
if $\{\mathrm{rank}(a_i)\mid i\in\omega\}$ is bounded by $m$
and $\{\mathrm{rank}(b_i)\mid i\in\omega\}$ by $n,$
then $\{\mathrm{rank}(a_i+b_i)\mid i\in\omega\}$ is bounded by $m+n.)$
It is easy to verify that for any nonprincipal
ultrafilter $\U$ on $\omega,$ the image of
this ideal in $\prod_i M_i(k)\,/\,\U$ is a proper nonzero ideal, so
unlike the $M_i(k),$ the ultraproduct ring is not simple.
(One can show that the ideals of this ring form an
uncountable chain.)
\vspace{.5em}

This finishes the material required for the preceding sections
of this paper.
We end by showing, for completeness's sake, how part of the
proof of Theorem~\ref{T.bigger} above can be strengthened.

In that proof, $K$ was an ultrapower of a field $k$ with respect to a
non-$\!\card(k)^+\!$-complete ultrafilter, and
we showed the degree $[K:k]$ to be uncountable by different
methods depending on whether $k$ was countable or uncountable:
in the former case by showing that $K$ had uncountable cardinality;
in the latter, by showing that it was transcendental over $k.$
Note that in the former situation, it follows that $K$ has transcendence
degree over $k$ equal to its cardinality; but our argument in
the latter case only proved transcendence degree $\geq 1.$
Can we similarly show in the second case that $k^I/\,\U$ has
transcendence degree over $k$ equal to its cardinality?

This is immediate if $\card(k^I/\,\U)>\card(k).$
To see this, note that it is easy to verify that whenever
$F$ is a transcendental extension of a field $E,$ one has
\begin{equation}\begin{minipage}[c]{35pc}\label{d.trdeg}
$\card(F)\ =\ \sup(\aleph_0,\ \card(E),\ \mathrm{tr.deg}_E(F)).$
\end{minipage}\end{equation}
Hence,
\begin{equation}\begin{minipage}[c]{35pc}\label{d.trdeg=}
If $\card(F)\ >\ \sup(\aleph_0,\ \card(E)),$ \ then
\ $\mathrm{tr.deg}_E(F)\ =\ \card(F).$
\end{minipage}\end{equation}
So if $\card(k^I/\,\U)>\card(k)\geq\aleph_0,$ we indeed get
$\mathrm{tr.deg}_k(k^I/\U)\ =\ \card(k^I/\U).$

Can the contrary case,
\begin{equation}\begin{minipage}[c]{35pc}\label{d.kI/U=k}
$\card(k^I/\,\U)\ =\ \card(k)$
\end{minipage}\end{equation}
occur for a non-$\!\card(k)^+\!$-complete ultrafilter $\U$?
Yes.
For instance, if we choose $k$ so that
$\card(k)=\linebreak[1]2^{\card(I)},$ then
$\card(k^I)=(2^{\card(I)})^{\card(I)}=2^{\card(I)}=\card(k),$ so
$\card(k^I/\,\U)$ certainly can't be larger.

We sketch below a proof that we nevertheless
have $\mathrm{tr.deg}_k(k^I/\,\U)=\card(k^I/\,\U),$
though under a slightly stronger hypothesis
that that of Theorem~\ref{T.bigger}; namely,
with the condition of that theorem that $\U$
be non-$\!\card(k)^+\!$-complete
strengthened to ``non-countably-complete''.
(This is strictly stronger only if
$\card(k)$ is $\geq$ some uncountable measurable cardinal.)

\begin{proposition}\label{P.tr_deg}
Let $k$ be an infinite field, and $\U$
a non-countably-complete ultrafilter on a set $I.$
Then $\mathrm{tr.deg}_k(k^I/\,\U)=\card(k^I/\,\U).$
\end{proposition}

\begin{proof}[Sketch of proof]
As noted above, the result is straightforward
unless~(\ref{d.kI/U=k}) holds, so assume~(\ref{d.kI/U=k}).

Let $k_0$ be any countable (possibly finite) subfield
of $k$ (e.g., its prime subfield), and note that since,
by~(\ref{d.kI/U=k}) and Theorem~\ref{T.bigger}, $k$ is
uncountable,~(\ref{d.trdeg=}) shows that
$\mathrm{tr.deg}_{k_0}(k)=\card(k).$
Let $\{s_{(\alpha)}\mid \alpha\in\card(k)\}$ be a
transcendence basis for $k$ over $k_0.$

Since $\U$ is not countably complete, we can
decompose $I$ into a countable family of disjoint $\!\U\!$-small
subsets $I_n$ $(n\in\omega).$
For each $\alpha\in\card(k),$
let $t_{(\alpha)}$ be the image in $k^I/\,\U$ of the
element of $k^I$ which on each $I_n$ has the constant
value $s_{(\alpha)}^n.$
We claim that these $t_{(\alpha)}$ are algebraically
independent over $k.$
Briefly, if they satisfied an algebraic dependence relation
over $k,$  they would satisfy such a relation over the
pure transcendental
subfield $k_0(s_{(\alpha)})_{\alpha\in\card(k)}.$
Hence, clearing denominators, we would get a polynomial
relation among the $t_{(\alpha)}$ with coefficients in
the polynomial ring $k_0[s_{(\alpha)}]_{\alpha\in\card(k)}.$
If this holds in $k^I/\,\U,$ then
the corresponding equation among elements of $k^I$
must hold at points of infinitely many $I_n,$ and by choosing
$n$ larger than the degrees in
the $s_{(\alpha)}$ of all the coefficients of the
given polynomial, one gets a contradiction.

This shows that the transcendence degree of
$k^I/\U$ over $k$ is at least
$\card(k),$ which, by~(\ref{d.kI/U=k}), equals $\card(k^I/\,\U).$
\end{proof}

We do not know whether in Proposition~\ref{P.tr_deg} the
assumption that $\U$ is not countably complete can be weakened
to ``not $\!\card(k)^+\!$-complete''.

\end{document}